\documentclass[a4paper,12pt, reqno]{amsart}
\usepackage{latexsym}
\usepackage{amssymb,amsfonts,amsmath,mathrsfs}
\begin{document}
\title{On the mean values of Dirichlet $L$-functions}
\author{H. M. BUI, J. P. KEATING}
\address{School of Mathematics, University of Bristol, Bristol, BS8 1TW}
\email{hm.bui@bristol.ac.uk, j.p.keating@bristol.ac.uk}

\begin{abstract}
We study the $2k$-th power moment of Dirichlet $L$-functions $L(s,\chi)$ at the centre of the critical strip $(s = 1/2)$, where the average is over all primitive characters $\chi\ (\textrm{mod}\ q)$. We extend to this case the hybrid Euler-Hadamard product results of Gonek, Hughes and Keating for the Riemann zeta-function. This allows us to recover conjectures for the moments based on random matrix models, incorporating the arithmetical terms in a
natural way.
\end{abstract}
\maketitle

\section{Introduction}

Let $L(s,\chi)$ denote a Dirichlet $L$-function. We shall here be interested in the $2k$-th power moment of $L(s,\chi)$ at the centre of the critical strip $(s=\frac{1}{2})$:
\begin{equation*}
\frac{1}{\varphi^{*}(q)}\sum_{\chi\ (\textrm{mod}\ q)}{\!\!\!\!\!\!\!}^{\textstyle{*}}\ |L({\scriptstyle{\frac{1}{2}}},\chi) |^{2k},
\end{equation*}
as $q\rightarrow\infty$, where $\sum^{*}$ denotes summation over all primitive characters $\chi$(mod $q$), and $\varphi^{*}(q)$ is the number of primitive characters. This is the $q$-analogue of the $2k$-th power moment of the Riemann zeta-function on the critical line.

It has long been known [\textbf{\ref{P}}] that
\begin{equation}\label{1}
\frac{1}{\varphi^{*}(q)}\sum_{\chi\ (\textrm{mod}\ q)}{\!\!\!\!\!\!\!}^{\textstyle{*}}\ |L({\scriptstyle{\frac{1}{2}}},\chi)|^{2}\sim\frac{\varphi(q)}{q}\log q,
\end{equation}
as $q\rightarrow \infty$. In 1981, Heath-Brown [\textbf{\ref{H-B2}}] proved that as $q\rightarrow\infty$,
\begin{equation}\label{2}
\frac{1}{\varphi^{*}(q)}\sum_{\chi\ (\textrm{mod}\
q)}{\!\!\!\!\!\!\!}^{\textstyle{*}}\ |L({\scriptstyle{\frac{1}{2}}},\chi)|^{4}=\frac{1}{2\pi^{2}}\prod_{p|q}\frac{(1-p^{-1})^3}{(1+p^{-1})}(\log q)^{4}
+O\bigg(\frac{2^{w(q)}q(\log q)^{3}}{\varphi^{*}(q)}\bigg),
\end{equation}
where $w(q)$ is the number of distinct prime factors of $q$. When $q$ has not too many prime factors, the error term in \eqref{2} is dominated by the main term, and this gives an asymptotic formula for the fourth moment. Recently, Soundararajan [\textbf{\ref{S}}] has improved Heath-Brown's result by showing that
\begin{equation}\label{3}
\frac{1}{\varphi^{*}(q)}\sum_{\chi\ (\textrm{mod}\ q)}{\!\!\!\!\!\!\!}^{\textstyle{*}}\ |L({\scriptstyle{\frac{1}{2}}},\chi)|^{4}\sim\frac{1}{2\pi^{2}}\prod_{p|q}\frac{(1-p^{-1})^3}{(1+p^{-1})}(\log q)^{4},
\end{equation}
as $q\rightarrow \infty$.

No asymptotic results for the $2k$-th moment have been proved when $k > 2$, though it has been conjectured [\textbf{\ref{KS1}},\textbf{\ref{CF}},\textbf{\ref{KS2}},\textbf{\ref{CFKRS}}] that the following holds.

\newtheorem{conj}{Conjecture}\begin{conj}
For $k$ fixed with $\emph{Re}k\geq0$,
\begin{equation*}
\frac{1}{\varphi^{*}(q)}\sum_{\chi\ (\emph{mod}\ q)}{\!\!\!\!\!\!\!}^{\textstyle{*}}\ |L({\scriptstyle{\frac{1}{2}}},\chi)|^{2k}\sim a(k)\frac{G^2(k+1)}{G(2k+1)}\prod_{p|q}\bigg(\sum_{m\geq0}\frac{d_k(p^m)^2}{p^m}\bigg)^{-1}(\log q)^{k^2},
\end{equation*}
as $q\rightarrow\infty$, where $G(z)$ is Barnes' $G$-function, $d_k(n)$ is the number of ways to represent $n$ as a product of $k$ factors, and
\begin{equation*}
a(k)=\prod_{p}\bigg(\bigg(1-\frac{1}{p}\bigg)^{k^2}\sum_{m\geq0}\frac{d_k(p^m)^2}{p^m}\bigg).
\end{equation*}
\end{conj}

\textbf{Remark}.\quad For $k\notin\mathbb{N}$, $d_k(p^n)$ is defined by $d_k(p^n)=\Gamma(n+k)/(n!\Gamma(k))$.\\

This conjecture was originally motivated by the random matrix model introduced by Keating and Snaith [\textbf{\ref{KS1}}], in which statistical properties of the $L$-functions are related to those of the characteristic polynomials of large random matrices. Specifically, let $U$ be an $N\times N$ unitary matrix. Denoting the eigenvalues of $U$ by $e^{i\theta_n}$, we see that the characteristic polynomial of $U$ is
\begin{equation*}
Z_N(U,\theta)=\prod_{n=1}^{N}(1-e^{i(\theta_n-\theta)}).
\end{equation*}
It was proved in [\textbf{\ref{KS1}}] that as $N\rightarrow\infty$,
\begin{equation}\label{4}
\mathbb{E}_N\big[|Z_N(U,\theta)|^{2k}\big]\sim\frac{G^2(k+1)}{G(2k+1)}N^{k^2},
\end{equation}
where the expectation value is computed with respect to Haar measure on $U(N)$. Equating the mean density of the eigenphases $\theta_n$ to the mean density of the $L$-function zeros corresponds to the identification $N\sim\log q$, and hence $N^{k^2}$ gives the right order for the $2k$-th moment of $L$-functions. However, the drawback of the model is the absence in \eqref{4} of the arithmetical factors $a(k)$ and the product over primes $p|q$ in Conjecture 1. These can be obtained straightforwardly from number-theoretical considerations [\textbf{\ref{CF}},\textbf{\ref{CFKRS}}], but then the random matrix contribution appears mysteriously. The question is how to treat the arithmetical and random matrix aspects on an equal footing.

Recently, it was shown by Gonek, Hughes and Keating [\textbf{\ref{GHK}}] that, using a smoothed form of the
explicit formula of Bombieri and Hejhal [\textbf{\ref{BH}}], one can approximate the Riemann zeta function
at a height $t$ on the critical line as a partial Euler product multiplied by a partial Hadamard product over the non-trivial zeros close to $1/2+it$. This suggests a statistical model for the zeta function in which the primes are involved in a natural way. The value distribution of the product over zeros is expected to be modelled by the characteristic polynomial of a large random unitary matrix, because it involves only local information about the zeros. Conjecturing the moments of this product using random matrix theory, calculating the moments of the product over the primes rigorously and making some assumptions (which can be proved in certain particular cases) about how the two products behave, Gonek, Hughes and Keating then reproduced the
conjecture about the $2k$-th moment of the zeta function first put forward by Keating and Snaith in [\textbf{\ref{KS1}}].

In this paper we show that the model introduced in [\textbf{\ref{GHK}}] can be adapted to the Dirichlet $L$-functions with primitive characters. Specifically, we mimic the results in [\textbf{\ref{GHK}}] to give the
following for the Dirichlet $L$-functions.

\newtheorem{theo}{Theorem}\begin{theo}
Let $u(x)$ be a real, non-negative, $C^\infty$-function with mass $1$ and compact support on $[e^{1-1/X},e]$. Set
\begin{equation*}
U(z)=\int_{0}^{\infty}u(x)E_{1}(z\log x)dx,
\end{equation*}
where $E_{1}(z)=\int_{z}^{\infty}e^{-w}/wdw$. Let $X\geq2$ be a real parameter. Then for $\chi$ a primitive character,
\begin{equation}\label{5}
L({\scriptstyle{\frac{1}{2}}},\chi)=P_{X}(\chi)Z_{X}(\chi)(1+O((\log X)^{-2})),
\end{equation}
where
\begin{equation}\label{6}
P_{X}(\chi)=\exp\bigg( \sum_{n\leq X}\frac{\Lambda(n)\chi(n)}{n^{1/2}\log n}\bigg),
\end{equation}
$\Lambda(n)$ is von Mangoldt's function, and
\begin{equation}\label{7}
Z_{X}(\chi)=\exp\bigg(-\sum_{\rho}U(({\scriptstyle{\frac{1}{2}}}-\rho)\log X)\bigg),
\end{equation}
where the sum is over the non-trivial zeros $\rho$ of $L(s,\chi)$.
\end{theo}

\textbf{Remark 1}.\quad The constant implied in the $O$-term is absolute. We can easily modify our result to handle $L(s,\chi)$ for all $\textrm{Re}s\geq0$ but then the constant in the error term is no longer absolute (cf. [\textbf{\ref{GHK}}, Theorem 1]).\\

\textbf{Remark 2}.\quad As discussed in [\textbf{\ref{GHK}}], $P_X(\chi)$ can be thought of as the Euler product for $L(\frac{1}{2},\chi)$ truncated to include primes $\lesssim X$, and $Z_X(\chi)$ can be thought of as the Hadamard product for $L(\frac{1}{2},\chi)$ truncated to include zeros within a distance $\lesssim 1/\log X$ of $s = \frac{1}{2}$. The parameter $X$ thus controls the relative contributions of the Euler and Hadamard products.\\

The proof of Theorem 1 is given in Section 2.

In Section 3, we evaluate the moments of $P_X(\chi)$ and prove the following.

\begin{theo}
Let $\delta>0$ and $k\geq0$, fixed. Suppose that $X,q\rightarrow\infty$ with $X\ll(\log q)^{2-\delta}$. Then
\begin{equation*}
\frac{1}{\varphi^{*}(q)}\sum_{\chi\ (\emph{mod}\ q)}{\!\!\!\!\!\!\!}^{\textstyle{*}}\ |P_{X}(\chi)|^{2k}
\sim a_k\prod_{\substack{p\leq X\\p|q}}\bigg(\sum_{m\geq 0}\frac{d_{k}(p^m)^2}{p^m}\bigg)^{-1}(e^{\gamma}\log X)^{k^2},
\end{equation*}
If we also have $w(q)\ll X^{1+\epsilon}$, where $w(q)$ is the number of distinct prime factors of $q$, then
\begin{equation*}
\frac{1}{\varphi^{*}(q)}\sum_{\chi\ (\emph{mod}\ q)}{\!\!\!\!\!\!\!}^{\textstyle{*}}\ |P_{X}(\chi)|^{2k}
\sim a_k\prod_{p|q}\bigg(\sum_{m\geq 0}\frac{d_{k}(p^m)^2}{p^m}\bigg)^{-1}(e^{\gamma}\log X)^{k^2}.
\end{equation*}
\end{theo}

The values of the moments of $Z_X(\chi)$ can be conjectured using random matrix theory as in [\textbf{\ref{GHK}}] (cf. Theorem 4 there). Instead of using matrices of size $N\sim\log T$, we would here average over all $N\times N$ unitary matrices with $N\sim\log q$ and so would have the following.

\begin{conj}
Let $k>-\frac{1}{2}$ be any real number. Suppose that $X,q\rightarrow\infty$ with $X\ll(\log q)^{2-\epsilon}$ and $w(q)\ll X^{1+\epsilon}$. Then
\begin{equation*}
\frac{1}{\varphi^{*}(q)}\sum_{\chi\ (\emph{mod}\ q)}{\!\!\!\!\!\!\!}^{\textstyle{*}}\ |Z_{X}(\chi)|^{2k}
\sim\frac{G^{2}(k+1)}{G(2k+1)}\bigg(\frac{\log q}{e^{\gamma}\log X}\bigg)^{k^2}.
\end{equation*}
\end{conj}

We note from Theorem 1 that
\begin{equation*}
L({\scriptstyle{\frac{1}{2}}},\chi)P_X(\chi)^{-1}=Z_X(\chi)(1+o(1)).
\end{equation*}
This allows us (in Section 5) to prove Conjecture 2 when $k=1$.

\begin{theo}
Let $\delta>0$. For $X,q\rightarrow\infty$ with $X\ll(\log q)^{2-\delta}$, we have
\begin{equation*}
\frac{1}{\varphi^{*}(q)}\sum_{\chi\ (\emph{mod}\ q)}{\!\!\!\!\!\!\!}^{\textstyle{*}}\ |Z_{X}(\chi)|^{2}\sim\prod_{\substack{p>X\\p|q}}\bigg(1-\frac{1}{p}\bigg)\frac{\log q}{e^{\gamma}\log X}.
\end{equation*}
If we also have $w(q)\ll X^{1+\epsilon}$, then
\begin{equation*}
\frac{1}{\varphi^{*}(q)}\sum_{\chi\ (\emph{mod}\ q)}{\!\!\!\!\!\!\!}^{\textstyle{*}}\ |Z_{X}(\chi)|^{2}\sim\frac{\log q}{e^{\gamma}\log X}.
\end{equation*}
\end{theo}

In the last section, we prove that Conjecture 2 also holds when $k = 2$, at least for the range $X\ll(\log\log q)^{2-\delta}$.

\begin{theo}
Let $\delta>0$. For $X,q\rightarrow\infty$ with $X\ll(\log\log q)^{2-\delta}$, we have
\begin{equation*}
\frac{1}{\varphi^{*}(q)}\sum_{\chi\ (\emph{mod}\ q)}{\!\!\!\!\!\!\!}^{\textstyle{*}}\ |Z_{X}(\chi)|^{4}
\sim\frac{1}{12}\prod_{\substack{p>X\\p|q}}\frac{(1-1/p)^{3}}{1+1/p}\bigg(\frac{\log q}{e^{\gamma}\log X}\bigg)^4.
\end{equation*}
If we also have $w(q)\ll X^{1+\epsilon}$, then
\begin{equation*}
\frac{1}{\varphi^{*}(q)}\sum_{\chi\ (\emph{mod}\ q)}{\!\!\!\!\!\!\!}^{\textstyle{*}}\ |Z_{X}(\chi)|^{4}
\sim\frac{1}{12}\bigg(\frac{\log q}{e^{\gamma}\log X}\bigg)^4.
\end{equation*}
\end{theo}

We remark that the condition on $X$ in Theorem 4 differs from that in the previous theorems. We believe that this is only a technical limitation and that the theorem should hold for a much larger range of values of $X$.

Combining the formulae for the second $(k = 1)$ and fourth $(k = 2)$ moments, \eqref{1} and \eqref{3}, with Theorem 2, Theorem 3 and Theorem 4, we see that, at least for the cases $k = 1$ and $k = 2$, when $X$ is not too large relative to $q$, the $2k$-th moment of $L(\frac{1}{2},\chi)$ is asymptotic to the product of the moments of $P_X(\chi)$ and $Z_X(\chi)$. We remark that it is rather interesting that one appears to need the condition $w(q)\ll X^{1+\epsilon}$ in Conjecture 2 for the moments of $Z_X(\chi)$ to coincide exactly with those of the characteristic polynomials of random matrices (cf. Theorems 3 and 4), but that even if this condition is not satisfied, when $k = 1$ and $k = 2$ the arithmetic dependence of these moments on $X$ cancels that of the corresponding moments of $P_X(\chi)$ so that \eqref{1} and \eqref{3} still follow from the product. We believe that this is true in general.

\begin{conj}
Let $k\geq0$ be any real number. Suppose that $X$ and $q\rightarrow\infty$ with $X\ll(\log q)^{2-\epsilon}$. Then
\begin{equation*}
\frac{1}{\varphi^{*}(q)}\sum_{\chi\ (\emph{mod}\ q)}{\!\!\!\!\!\!\!}^{\textstyle{*}}\ |L({\scriptstyle{\frac{1}{2}}},\chi)|^{2k}\sim\bigg(\frac{1}{\varphi^{*}(q)}\sum_{\chi\ (\emph{mod}\ q)}{\!\!\!\!\!\!\!}^{\textstyle{*}}\ |P_{X}(\chi)|^{2k}\bigg)\bigg(\frac{1}{\varphi^{*}(q)}\sum_{\chi\ (\emph{mod}\ q)}{\!\!\!\!\!\!\!}^{\textstyle{*}}\ |Z_{X}(\chi)|^{2k}\bigg).
\end{equation*}
\end{conj}

This conjecture, together with Theorem 2 and the random matrix model for $Z_{X}(\chi)$, implies Conjecture 1.

\section{Proof of Theorem 1}

Similar to [\textbf{\ref{GHK}}] (cf. Section 2 there), we have the following result.

\newtheorem{lemm}{Lemma}\begin{lemm}
Let $u(x)$ be a real, non-negative, $C^{\infty}$ function with mass $1$ and compact support on $[e^{1-1/X},e]$. Let $v(t)=\int_{t}^{\infty}u(x)dx$ and let $\tilde{u}$ be the Mellin transform of $u$. Then for $\chi$ a primitive character and $s$ not a zero of $L(s,\chi)$, we have
\begin{eqnarray*}
&&-\frac{L'}{L}(s,\chi)=\sum_{n=2}^{\infty}\frac{\Lambda(n)\chi(n)}{n^s}v(e^{\log n/\log X})-\sum_{\rho}\frac{\tilde{u}(1-(s-\rho)\log X)}{s-\rho}\nonumber\\
&&\qquad\qquad\qquad\qquad-\sum_{m=0}^{\infty}\frac{\tilde{u}(1-(s+\mathfrak{a}+2m)\log X)}{s+\mathfrak{a}+2m},
\end{eqnarray*}
where $\mathfrak{a}$ is defined by $\chi(-1)=(-1)^{\mathfrak{a}}$ and the sum over $\rho$ runs over all the non-trivial zeros of $L(s,\chi)$.
\end{lemm}

This lemma can be proved in a familiar way [\textbf{\ref{BH}}], beginning with the integral
\begin{equation*}
\frac{1}{2\pi i}\int_{(c)}\frac{L'}{L}(z+s)\tilde{u}(1+z\log X)\frac{dz}{z},
\end{equation*}
where $c=\max\{2,2-\textrm{Re}s\}$.

Following the arguments in [\textbf{\ref{GHK}}], we can integrate the formula in Lemma 1 to give a formula
for $L(s,\chi)$: we have, for $s$ not equal to one of the zeros of the Dirichlet $L$-function and $\textrm{Re}s\geq0$,

\begin{eqnarray}\label{8}
L(s,\chi)&=&\exp\bigg(\sum_{n=2}^{\infty}\frac{\Lambda(n)\chi(n)}{n^s\log n}v(e^{\log n/\log X})\bigg)Z_{X}(\chi)\nonumber\\
&&\qquad\qquad\times\exp\bigg(-\sum_{m=0}^{\infty}U((s+\mathfrak{a}+2m)\log X)\bigg).
\end{eqnarray}
To remove the former restriction on $s$, we note that we may interpret $\exp(-U(z))$ to be asymptotic to $Cz$ for some constant $C$ as $z\rightarrow0$, so both sides of \eqref{8} vanish at the zeros. Thus \eqref{8} holds for all $\textrm{Re}s\geq0$. Let $s=1/2$, and observe that
\begin{eqnarray*}
U(({\scriptstyle{\frac{1}{2}}}+\mathfrak{a}+2m)\log X)&=&\int_{e^{1-1/X}}^{e}u(y)E_{1}(({\scriptstyle{\frac{1}{2}}}+\mathfrak{a}+2m)\log X\log y)dy\nonumber\\
&\ll&\frac{1}{(1/2+\mathfrak{a}+2m)^{2}(\log X)^2},
\end{eqnarray*}
as for $x>0$, $E_{1}(x)=\int_{x}^{\infty}e^{-w}/wdw\ll x^{-2}$. Hence
\begin{equation}\label{9}
\exp\bigg(-\sum_{m=0}^{\infty}U(({\scriptstyle{\frac{1}{2}}}+\mathfrak{a}+2m)\log X)\bigg)=1+O((\log X)^{-2}).
\end{equation}
Also, since $v(e^{\log n/\log X})=1$ for $n\leq X^{1-1/X}$, the first factor in \eqref{8} is
\begin{eqnarray}\label{10}
&&P_{X}(\chi)\exp\bigg(\sum_{X^{1-1/X}\leq n\leq X}\frac{\Lambda(n)\chi(n)}{n^{1/2}\log n}(v(e^{\log n/\log X})-1)\bigg)\nonumber\\
&&\qquad\qquad\qquad=P_{X}(\chi)\exp\bigg(O\bigg(\sum_{X^{1-1/X}\leq n\leq X}\frac{1}{n^{1/2}}\bigg)\bigg)\nonumber\\
&&\qquad\qquad\qquad=P_{X}(\chi)\exp(O(X^{-1/2}\log X))\nonumber\\
&&\qquad\qquad\qquad=P_{X}(\chi)(1+O(X^{-1/2}\log X)).
\end{eqnarray}
Theorem 1 then follows from \eqref{8}, \eqref{9} and \eqref{10}.

\section{Proof of Theorem 2}

Our strategy is to express $P_X(\chi)^k$ as a Dirichlet polynomial and use the orthogonality relation of Dirichlet characters. We require some lemmas.

\begin{lemm}
Let $k\geq0$ be fixed. Suppose that $X,q\rightarrow\infty$ with $w(q)\ll X^{1+\epsilon}$. Then
\begin{equation*}
\prod_{\substack{p>X\\p|q}}\bigg(\sum_{m\geq 0}\frac{d_{k}(p^m)^2}{p^m}\bigg)=1+o(1).
\end{equation*}
\end{lemm}
\begin{proof}
Let $l=w(q)$ and $p_{1}<p_{2}<\ldots<p_{l}$ be the first $l$ primes after $X$. By the prime number theorem, we have $p_l/\log p_l\sim l+X/\log X$. So
\begin{equation*}
\log p_{l}\leq(1+o(1))\log\bigg(l+\frac{X}{\log X}\bigg)\leq(1+o(1))\log X.
\end{equation*}
Thus
\begin{eqnarray*}
\sum_{\substack{p>X\\p|q}}\frac{1}{p}\leq\sum_{j\leq l}\frac{1}{p_{j}}&=&\log\log p_{l}-\log\log p_{1}+O\bigg(\frac{1}{\log X}\bigg)\nonumber\\
&=&\log\log p_{l}-\log\log X+O\bigg(\frac{1}{\log X}\bigg)=o(1).
\end{eqnarray*}
The result easily follows using the bound $d_{k}(p^m)=O(p^{\epsilon})$.
\end{proof}

\begin{lemm}
Let
\begin{equation*}
P_{X}(s,\chi)=\exp\bigg( \sum_{n\leq X}\frac{\Lambda(n)\chi(n)}{n^{s}\log n}\bigg),
\end{equation*}
so $P_{X}(\chi)=P_{X}({\scriptstyle{\frac{1}{2}}},\chi)$, and let $P_{X}^{*}(\chi)=P_{X}^{*}({\scriptstyle{\frac{1}{2}}},\chi)$, where
\begin{equation*}
P_{X}^{*}(s,\chi)=\prod_{p\leq X}\bigg( 1-\frac{\chi(p)}{p^s}\bigg) ^{-1}\prod_{\sqrt{X}<p\leq X}\bigg( 1+\frac{\chi(p)^2}{2p^{2s}}\bigg) ^{-1}.
\end{equation*}
Then for any $k\in\mathbb{R}$ we have
\begin{equation*}
P_{X}(s,\chi)^{k}=P_{X}^{*}(s,\chi)^{k}\bigg(1+O_{k}\bigg(\frac{1}{\log X}\bigg)\bigg),
\end{equation*}
uniformly for $\sigma\geq1/2$.
\end{lemm}
\begin{proof}
Let $N_{p}=[\log X/\log p]$, the integer part of $\log X/\log p$. We have
\begin{equation*}
P_{X}(s,\chi)^{k}=\exp\bigg(k\sum_{p\leq X}\sum_{1\leq j\leq N_{p}}\frac{\chi(p)^j}{jp^{js}}\bigg),
\end{equation*}
and
\begin{equation*}
P_{X}^{*}(s,\chi)^{k}=\exp\bigg(k\sum_{p\leq X}\sum_{j\geq1}\frac{\chi(p)^j}{jp^{js}}+k\sum_{\sqrt{X}<p\leq X}\sum_{j\geq1}\frac{(-1)^{j}\chi(p)^{2j}}{j2^{j}p^{2js}}\bigg).
\end{equation*}
Therefore
\begin{equation*}
P_{X}(s,\chi)^{k}P_{X}^{*}(s,\chi)^{-k}=\exp\bigg(-k\sum_{p\leq X}\sum_{j>N_{p}}\frac{\chi(p)^j}{jp^{js}}-k\sum_{\sqrt{X}<p\leq X}\sum_{j\geq1}\frac{(-1)^{j}\chi(p)^{2j}}{j2^{j}p^{2js}}\bigg).
\end{equation*}
We note that $N_{p}=1$ for $\sqrt{X}<p\leq X$, and the $j=2$ term for these primes in the first double sum cancels the $j=1$ term in the second. Thus the expression in the exponent is
\begin{eqnarray*}
&\ll&|k|\bigg(\sum_{p\leq\sqrt{X}}\frac{1}{p^{\sigma(N_{p}+1)}}+\sum_{\sqrt{X}<p\leq X}\frac{1}{p^{3\sigma}}\bigg)\nonumber\\
&\ll&|k|\bigg(X^{-\sigma}\sum_{p\leq\sqrt{X}}1+\sum_{\sqrt{X}<p\leq X}\frac{1}{p^{3/2}}\bigg)\\
&\ll&|k|\bigg(\frac{1}{\log X}+\frac{1}{X^{1/4}\log X}\bigg)\ll\frac{|k|}{\log X}.
\end{eqnarray*}
Hence $P_{X}(s,\chi)^{k}P_{X}^{*}(s,\chi)^{-k}=1+O_{k}(1/\log X)$ as required.
\end{proof}

The next lemma is standard (see [\textbf{\ref{S}}, Lemma 1]).

\begin{lemm}
For $(mn,q)=1$, we have
\begin{equation*}
\sum_{\chi\ (\emph{mod}\ q)}{\!\!\!\!\!\!\!}^{\textstyle{*}}\ \ \chi(m)\overline\chi(n)=\sum_{\substack{h|q\\h|(m-n)}}\varphi(h)\mu(q/h),
\end{equation*}
and if we restrict to characters of a given sign $\mathfrak{a}$, then
\begin{equation*}
\sum_{\substack{\chi\ (\emph{mod}\ q)\\\chi(-1)=(-1)^{\mathfrak{a}}}}{\!\!\!\!\!\!\!\!\!\!}^{\textstyle{*}}\ \ \chi(m)\overline\chi(n)=\frac{1}{2}\sum_{\substack{h|q\\h|(m-n)}}\varphi(h)\mu(q/h)+\frac{(-1)^{\mathfrak{a}}}{2}\sum_{\substack{h|q\\h|(m+n)}}\varphi(h)\mu(q/h).
\end{equation*}
\end{lemm}

\textit{Remark}.\quad In particular, for $m=n$ we obtain the formula for $\varphi^{*}(q)$:
\begin{equation*}
\varphi^{*}(q)=q\prod_{p||q}\bigg(1-\frac{2}{p}\bigg)\prod_{p^2|q}\bigg(1-\frac{1}{p}\bigg)^2.
\end{equation*}

We now proceed with the proof of Theorem 2. We write $P_{X}^{*}(s,\chi)^{k}$ as a Dirichlet series
\begin{equation}\label{11}
\sum_{n=1}^{\infty}\frac{\alpha_{k}(n)\chi(n)}{n^s}=\prod_{p\leq X}\bigg( 1-\frac{\chi(p)}{p^s}\bigg) ^{-k}\prod_{\sqrt{X}<p\leq X}\bigg( 1+\frac{\chi(p)^2}{2p^{2s}}\bigg) ^{-k}.
\end{equation}
We note that $\alpha_{k}(n)\in\mathbb{R}$, and if we denote by $S(X)$ the set of $X$-smooth numbers, that is
\begin{equation*}
S(X)=\{n\in\mathbb{N}:p|n\rightarrow p\leq X\},
\end{equation*}
then $\alpha_{k}(n)$ is multiplicative, and $\alpha_{k}(n)=0$ if $n\notin S(X)$. We also have $\alpha_{k}(n)=d_{k}(n)$ if $n\in S(\sqrt{X})$, and $\alpha_{k}(p)=d_{k}(p)=k$. Moreover, by comparing $(1-\chi(p)p^{-s})^{-k}$ with $(1-\chi(p)p^{-s})^{-k}(1+\chi(p)^2p^{-2s}/2)^{-k}$, we easily find that $|\alpha_{k}(n)|\leq d_{3|k|/2}(n)$.

We now truncate the series, for $s=1/2$, at $q^{\theta}$, where $\theta>0$ will be chosen later. We have
\begin{equation*}
\sum_{n\in S(X)}\frac{\alpha_{k}(n)\chi(n)}{\sqrt{n}}=\sum_{\substack{n\in S(X)\\n\leq  q^{\theta}}}\frac{\alpha_{k}(n)\chi(n)}{\sqrt{n}}+O\bigg(\sum_{\substack{n\in S(X)\\n> q^{\theta}}}\frac{d_{3|k|/2}(n)}{\sqrt{n}}\bigg).
\end{equation*}
The $O$-term is
\begin{eqnarray*}
&\ll&\sum_{n\in S(X)}\bigg(\frac{n}{q^{\theta}}\bigg)^{\delta/4}\frac{d_{3|k|/2}(n)}{\sqrt{n}}=q^{-\delta\theta/4}\prod_{p\leq X}\bigg(1-\frac{1}{p^{1/2-\delta/4}}\bigg)^{-3|k|/2}\nonumber\\
&\ll&q^{-\delta\theta/4}\exp\bigg(O_{k}\bigg(\frac{X^{1/2+\delta/4}}{\log X}\bigg)\bigg)\ll q^{-\delta\theta/4}\exp\bigg(O_{k}\bigg(\frac{\log q}{\log\log q}\bigg)\bigg)\ll q^{-\delta\theta/5},
\end{eqnarray*}
since 
\begin{equation*}
X^{1/2+\delta/4}\ll(\log q)^{(2-\delta)(1/2+\delta/4)}\ll\log q.
\end{equation*}
So Lemma 3 yields
\begin{equation}\label{12}
P_{X}(\chi)^{k}=\bigg(\sum_{\substack{n\in S(X)\\n\leq  q^{\theta}}}\frac{\alpha_{k}(n)\chi(n)}{\sqrt{n}}+O(q^{-\delta\theta/5})\bigg)\bigg(1+O_{k}\bigg(\frac{1}{\log X}\bigg)\bigg).
\end{equation}
From Lemma 4, we have
\begin{eqnarray*}
\frac{1}{\varphi^{*}(q)}\sum_{\chi\ (\textrm{mod}\ q)}{\!\!\!\!\!\!\!}^{\textstyle{*}}\ \ \bigg|\sum_{\substack{n\in S(X)\\n\leq q^{\theta}}}\frac{\alpha_{k}(n)\chi(n)}{\sqrt{n}}\bigg|^{2}&=&\frac{1}{\varphi^{*}(q)}\sum_{\substack{mn\in S(X)\\m,n\leq q^{\theta}}}\frac{\alpha_{k}(m)\alpha_{k}(n)}{\sqrt{mn}}\sum_{\chi\ (\textrm{mod}\ q)}{\!\!\!\!\!\!\!}^{\textstyle{*}}\ \chi(m)\overline\chi(n)\nonumber\\
&=&\frac{1}{\varphi^{*}(q)}\sum_{h|q}\varphi(h)\mu(q/h)\sum_{\substack{mn\in S(X)\\m,n\leq q^{\theta}\\h|(m-n)\\(mn,q)=1}}\frac{\alpha_{k}(m)\alpha_{k}(n)}{\sqrt{mn}}\nonumber\\
&=&T_{1}+T_{2},
\end{eqnarray*}
where $T_{1}$ consists of the diagonal terms $m=n$ and $T_{2}$ is the sum of the remaining terms. 

We first estimate $T_{1}$. We have
\begin{equation*}
T_{1}=\sum_{\substack{n\in S(X)\\n\leq q^{\theta}\\(n,q)=1}}\frac{\alpha_{k}(n)^2}{n}.
\end{equation*}
Using the method above, we may extend the sum to infinity with a gain of at most $O(q^{-\theta/3})$. So
\begin{equation*}
T_{1}=\sum_{\substack{n\in S(X)\\(n,q)=1}}\frac{\alpha_{k}(n)^2}{n}+O(q^{-\theta/3}).
\end{equation*}
We note again that $\alpha_{k}(n)=d_{k}(n)$ for $n\in S(\sqrt{X})$ or when $n$ is prime between $(\sqrt{X},X]$. Hence we can write the sum as
\begin{eqnarray*}
\prod_{\substack{p\leq X\\p\nmid q}}\bigg(\sum_{m\geq0}\frac{\alpha_{k}(p^m)^2}{p^m}\bigg)&=&\prod_{\substack{p\leq X\\p\nmid q}}\bigg(\sum_{m\geq0}\frac{d_{k}(p^m)^2}{p^m}\bigg)\prod_{\substack{\sqrt{X}<p\leq X\\p\nmid q}}\frac{1+\frac{d_{k}(p)^2}{p}+\sum_{m\geq2}\frac{\alpha_{k}(p^m)^2}{p^m}}{\sum_{m\geq0}\frac{d_{k}(p^m)^2}{p^m}}\\
&=&\prod_{\substack{p\leq X\\p\nmid q}}\bigg(\sum_{m\geq 0}\frac{d_{k}(p^m)^2}{p^m}\bigg)\prod_{\substack{\sqrt{X}<p\leq X\\p\nmid q}}\bigg(1+O_{k}\bigg(\frac{1}{p^2}\bigg)\bigg)\\
&=&(1+o(1))\prod_{p\leq X}\bigg(\sum_{m\geq 0}\frac{d_{k}(p^m)^2}{p^m}\bigg)\prod_{\substack{p\leq X\\p|q}}\bigg(\sum_{m\geq 0}\frac{d_{k}(p^m)^2}{p^m}\bigg)^{-1}.
\end{eqnarray*}
The first product, by Mertens' theorem, is equal to
\begin{eqnarray*}
&&\prod_{p\leq X}\bigg(\bigg(1-\frac{1}{p}\bigg)^{k^2}\sum_{m\geq0}\frac{d_{k}(p^m)^2}{p^m}\bigg)\prod_{p\leq X}\bigg(1-\frac{1}{p}\bigg)^{-k^2}\nonumber\\
&&\qquad\qquad=(1+o(1))\prod_{p>X}\bigg(\bigg(1-\frac{1}{p}\bigg)^{k^2}\sum_{m\geq0}\frac{d_{k}(p^m)^2}{p^m}\bigg)^{-1}a_k(e^{\gamma}\log X)^{k^2}\nonumber\\
&&\qquad\qquad=(1+o(1))\prod_{p>X}\bigg(1+O_{k}\bigg(\frac{1}{p^2}\bigg)\bigg)a_k(e^{\gamma}\log X)^{k^2}\nonumber\\
&&\qquad\qquad=(1+o(1))a_k(e^{\gamma}\log X)^{k^2}.
\end{eqnarray*}
Thus
\begin{equation}\label{13}
T_{1}=(1+o(1))a_k\prod_{\substack{p\leq X\\p|q}}\bigg(\sum_{m\geq 0}\frac{d_{k}(p^m)^2}{p^m}\bigg)^{-1}(e^{\gamma}\log X)^{k^2}.
\end{equation}
We note that $T_{1}\gg (\varphi(q)/q)^{k^2}(\log X)^{k^2}$. 

To estimate $T_{2}$, we note that for $m\equiv n(\textrm{mod}\ h)$, $m\ne n$ and $m,n\leq q^{\theta}$, we can restrict the sum over $h$ to $h\leq q^{\theta}$. So
\begin{eqnarray}\label{14}
T_{2}&\ll&\frac{1}{\varphi^{*}(q)}\sum_{\substack{h|q\\h\leq q^{\theta}}}\varphi(h)\bigg(\sum_{\substack{mn\in S(X)\\m,n\leq q^{\theta}}}\frac{|\alpha_{k}(m)\alpha_{k}(n)|}{\sqrt{mn}}\bigg)\nonumber\\
&\ll&\frac{q^{2\theta}}{\varphi^{*}(q)}\bigg(\sum_{n\in S(X)}\frac{d_{3|k|/2}(n)}{\sqrt{n}}\bigg)^2\ll \frac{q^{2\theta}}{\varphi^{*}(q)}\prod_{p\leq X}\bigg(1-\frac{1}{\sqrt{p}}\bigg)^{-3|k|}\nonumber\\
&\ll&\frac{q^{2\theta}}{\varphi^{*}(q)}e^{3|k|\sqrt{X}}\ll \frac{q^{3\theta}}{\varphi^{*}(q)}\ll q^{-1+3\theta}\bigg(\frac{q}{\varphi(q)}\bigg)^2.
\end{eqnarray}
The last line follows by the prime number theorem and because 
\begin{equation*}
3|k|\sqrt{X}\ll(\log q)^{1-\delta/2}\ll\theta\log q.
\end{equation*}

Since $q/\varphi(q)\ll\log\log q$, choosing $\theta=\frac{1}{4}$, we find that \eqref{14}, together with \eqref{13} and \eqref{12}, complete the proof of the first part of the theorem.

The second statement of the theorem then easily follows by Lemma 2.

\section{Functional equations}

For a primitive character $\chi\ (\textrm{mod}\ q)$, let $\mathfrak{a}$ be given by $\chi(-1)=(-1)^{\mathfrak{a}}$. Define
\begin{equation*}
\Lambda({\scriptstyle{\frac{1}{2}}}+s,\chi)=\bigg(\frac{q}{\pi}\bigg)^{s/2}\Gamma\bigg(\frac{s+1/2+\mathfrak{a}}{2}\bigg)L({\scriptstyle{\frac{1}{2}}}+s,\chi).
\end{equation*}
This is an entire function and it satisfies the functional equation
\begin{equation*}
\Lambda({\scriptstyle{\frac{1}{2}}}+s,\chi)=\frac{\tau(\chi)}{i^{\mathfrak{a}}\sqrt{q}}\Lambda({\scriptstyle{\frac{1}{2}}}-s,\overline\chi).
\end{equation*}
For $c>\frac{1}{2}$, we consider
\begin{equation*}
A(\chi):=\frac{1}{2\pi i}\int_{(1)}\frac{\Lambda({\scriptstyle{\frac{1}{2}}}+s,\chi)\Lambda({\scriptstyle{\frac{1}{2}}}+s,\overline\chi)}{\Gamma(\frac{1/2+\mathfrak{a}}{2})^2}\frac{ds}{s}.
\end{equation*}
Moving the line of integration to $\Re s=-c$, and applying Cauchy's theorem and the functional equation, we deduce that $A(\chi)=|L(\frac{1}{2},\chi)|^2-A(\chi)$. Also, expanding $L(\frac{1}{2}+s,\chi)L(\frac{1}{2}+s,\overline\chi)$ in a Dirichlet series and integrating termwise we get
\begin{equation*}
A(\chi)=\sum_{a,b\geq1}\frac{\chi(a)\overline{\chi}(b)}{\sqrt{ab}}W_{\mathfrak{a}}\bigg(\frac{\pi ab}{q}\bigg),
\end{equation*}
where
\begin{equation*}
W_{\mathfrak{a}}(x)=\frac{1}{2\pi i}\int_{(1)}\bigg(\frac{\Gamma(\frac{1/2+s+\mathfrak{a}}{2})}{\Gamma(\frac{1/2+\mathfrak{a}}{2})}\bigg)^{2}x^{-s}\frac{ds}{s}.
\end{equation*}
We have $W_{\mathfrak{a}}(x)=O_c(x^{-c})$ and also, by moving the line of integration to $c=-\frac{1}{2}+\epsilon$, we have $W_{\mathfrak{a}}(x)=1+O(x^{1/2+\epsilon})$.

Let $Z=q/2^{\omega(q)}$ and decompose $A(\chi)$ as $B(\chi)+C(\chi)$ where
\begin{equation*}
B(\chi)=\sum_{\substack{a,b\geq1\\ab\leq Z}}\frac{\chi(a)\overline{\chi}(b)}{\sqrt{ab}}W_{\mathfrak{a}}\bigg(\frac{\pi ab}{q}\bigg),
\end{equation*}
and
\begin{equation*}
C(\chi)=\sum_{\substack{a,b\geq1\\ab>Z}}\frac{\chi(a)\overline{\chi}(b)}{\sqrt{ab}}W_{\mathfrak{a}}\bigg(\frac{\pi ab}{q}\bigg).
\end{equation*}
Our aim in the next two sections is to evaluate the first and the second moments of $B(\chi)|P_X(\chi)^{-1}|^2$ and $C(\chi)|P_X(\chi)^{-1}|^2$. Theorems 3 and 4 will then easily follow.

\section{Proof of Theorem 3}

We begin with some lemmas.

\begin{lemm}
Let $q$ be a positive integer and $x\geq2$. Then
\begin{equation*}
\sum_{\substack{n\leq x\\(n,q)=1}}\frac{1}{n}=\frac{\varphi(q)}{q}(\log x+O(1+\log\omega(q)))+O\bigg(\frac{2^{\omega(q)}\log x}{x}\bigg).
\end{equation*}
\end{lemm}
\begin{proof}
We have
\begin{eqnarray*}
\sum_{\substack{n\leq x\\(n,q)=1}}\frac{1}{n}=\sum_{d|q}\mu(d)\sum_{\substack{n\leq x\\d|n}}\frac{1}{n}&=&\sum_{\substack{d\leq x\\d|q}}\frac{\mu(d)}{d}\bigg(\log\frac{x}{d}+\gamma+O\bigg(\frac{d}{x}\bigg)\bigg)\nonumber\\
&=&\sum_{d|q}\frac{\mu(d)}{d}\bigg(\log\frac{x}{d}+\gamma\bigg)+O\bigg(\frac{2^{w(q)}\log x}{x}\bigg)\\
&=&\frac{\varphi(q)}{q}\bigg(\log x+\gamma+\sum_{p|q}\frac{\log p}{p-1}\bigg)+O\bigg(\frac{2^{w(q)}\log x}{x}\bigg).
\end{eqnarray*}
Since $\sum_{p|q}\log p/(p-1)\ll 1+\log\omega(q)$, the lemma follows.
\end{proof}

\begin{lemm}
Let $m$, $n$, $h$ be three positive integers, $(mn,h)=1$, and $Z_{1}\geq 2$. Then
\begin{equation*}
E_{1}:=\sum_{\substack{Z_{1}\leq ab<2Z_{1}\\am\equiv\pm bn(\emph{mod}\ h)\\am\ne bn\\(ab,h)=1}}1\ll\frac{Z_{1}mn}{h}\log (Z_{1}mn).
\end{equation*}
\end{lemm}
\begin{proof}
It is clear that we only need to consider $am>bn$. Let $am=lh\pm bn$, where $l\geq1$. We have $l\leq(am+bn)/h\leq 2Z_{1}(m+n)/h$ and $b\leq2Z_{1}/a\leq 4Z_{1}m/lh$. So
\begin{equation*}
E_{1}\ll\sum_{1\leq l\leq\frac{2Z_{1}(m+n)}{h}}\sum_{\substack{b\leq\frac{4Z_{1}m}{lh}\\(b,h)=1}}1\ll\frac{Z_{1}m}{h}\sum_{1\leq l\leq\frac{2Z_{1}(m+n)}{h}}\frac{1}{l}\ll\frac{Z_{1}mn}{h}\log (Z_{1}mn).
\end{equation*}
The proof is complete.
\end{proof}

We next prove two propositions.

\newtheorem{prop}{Proposition}\begin{prop}
Let $\delta>0$. Suppose we have $X,q\rightarrow\infty$ with $X\ll(\log q)^{2-\delta}$. Then
\begin{equation*}
I=\frac{1}{\varphi^{*}(q)}\sum_{\chi\ (\emph{mod}\ q)}{\!\!\!\!\!\!\!}^{\textstyle{*}}\ B(\chi)\bigg|\sum_{\substack{n\in S(X)\\n\leq q^{1/10}}}\frac{\alpha_{-1}(n)\chi(n)}{\sqrt{n}}\bigg|^2\sim\frac{1}{2}\prod_{\substack{p>X\\p|q}}\bigg(1-\frac{1}{p}\bigg)\frac{\log q}{e^{\gamma}\log X}.
\end{equation*}
\end{prop}
\begin{proof}
We have
\begin{eqnarray*}
I&=&\frac{1}{\varphi^{*}(q)}\sum_{\substack{ab\leq Z\\mn\in S(X)\\m,n\leq q^{1/10}}}\frac{\alpha_{-1}(m)\alpha_{-1}(n)}{\sqrt{abmn}}\sum_{\mathfrak{a}\in\{0,1\}}W_{\mathfrak{a}}\bigg(\frac{\pi ab}{q}\bigg)\sum_{\substack{\chi\ (\textrm{mod}\ q)\\\chi(-1)=(-1)^{\mathfrak{a}}}}{\!\!\!\!\!\!\!\!\!}^{\textstyle{*}}\ \chi(am)\overline{\chi}(bn)\nonumber\\
&=&I_{1}+I_{2},
\end{eqnarray*}
where $I_{1}$ and $I_{2}$ are, respectively, the diagonal and the off-diagonal. We first consider $I_1$:
\begin{equation*}
I_{1}=\frac{1}{2}\sum_{\mathfrak{a}\in\{0,1\}}\sum_{\substack{ab\leq Z\\mn\in S(X)\\m,n\leq q^{1/10}\\am=bn\\(abmn,q)=1}}\frac{\alpha_{-1}(m)\alpha_{-1}(n)}{\sqrt{abmn}}W_{\mathfrak{a}}\bigg(\frac{\pi ab}{q}\bigg).
\end{equation*}
Since $am=bn$, we can write $m=ur$, $n=us$, $a=vs$ and $b=vr$, where $(r,s)=1$. Then
\begin{eqnarray}\label{27}
I_{1}&=&\frac{1}{2}\sum_{\mathfrak{a}\in\{0,1\}}\sum_{\substack{v^{2}rs\leq Z\\urs\in S(X)\\ur,us\leq q^{1/10}\\(uvrs,q)=1\\(r,s)=1}}\frac{\alpha_{-1}(ur)\alpha_{-1}(us)}{uvrs}W_{\mathfrak{a}}\bigg(\frac{\pi v^{2}rs}{q}\bigg)\nonumber\\
&=&\frac{1}{2}\sum_{\mathfrak{a}\in\{0,1\}}\sum_{\substack{urs\in S(X)\\ur,us\leq q^{1/10}\\(urs,q)=1\\(r,s)=1}}\frac{\alpha_{-1}(ur)\alpha_{-1}(us)}{urs}\sum_{\substack{v\leq\sqrt{\frac{Z}{rs}} \\(v,q)=1}}\frac{1}{v}W_{\mathfrak{a}}\bigg(\frac{\pi v^{2}rs}{q}\bigg)
\end{eqnarray}
The sum over $v$ is
\begin{equation*}
\sum_{\substack{v\leq\sqrt{\frac{Z}{rs}}\\(v,q)=1}}\frac{1}{v}\bigg(1+O\bigg(\frac{\sqrt{v}(rs)^{1/4}}{q^{1/4}}\bigg)\bigg)=\sum_{\substack{v\leq\sqrt{\frac{Z}{rs}}\\(v,q)=1}}\frac{1}{v}+O(2^{-\omega(q)/4}),
\end{equation*}
which is, by Lemma 5, equal to
\begin{equation*}
\frac{\varphi(q)}{2q}\log\frac{Z}{rs}+O(1+\log\omega(q))=\frac{\varphi(q)}{2q}((1+o(1))\log q-\log(rs)).
\end{equation*}
Hence the right-hand side in \eqref{27} breaks into, say, $I_{1}^{M}+I_{1}^{E}$, where
\begin{equation*}
I_{1}^{M}=(1+o(1))\frac{\varphi(q)}{2q}\log q\sum_{\substack{urs\in S(X)\\ur,us\leq q^{1/10}\\(urs,q)=1\\(r,s)=1}}\frac{\alpha_{-1}(ur)\alpha_{-1}(us)}{urs},
\end{equation*}
and
\begin{equation*}
I_{1}^{E}=-\frac{\varphi(q)}{2q}\sum_{\substack{urs\in S(X)\\ur,us\leq q^{1/10}\\(urs,q)=1\\(r,s)=1}}\frac{\alpha_{-1}(ur)\alpha_{-1}(us)\log (rs)}{urs}.
\end{equation*}
As in [\textbf{\ref{GHK}}] (cf. Section 5 there), we have $I_{1}^{E}\ll(\log X)^{10}$. Furthermore, since $\sum_{d|n}\varphi(d)=n$, we have
\begin{eqnarray*}
I_{1}^{M}&=&(1+o(1))\frac{\varphi(q)}{2q}\log q\sum_{\substack{mn\in S(X)\\m,n\leq q^{1/10}\\(mn,q)=1}}\frac{\alpha_{-1}(m)\alpha_{-1}(n)}{mn}(m,n)\nonumber\\
&=&(1+o(1))\frac{\varphi(q)}{2q}\log q\sum_{\substack{mn\in S(X)\\m,n\leq q^{1/10}\\(mn,q)=1}}\frac{\alpha_{-1}(m)\alpha_{-1}(n)}{mn}\sum_{\substack{u|m\\u|n}}\varphi(u)\nonumber\\
&=&\label{2.17}(1+o(1))\frac{\varphi(q)}{2q}\log q\sum_{\substack{u\in S(X)\\u\leq q^{1/10}\\(u,q)=1}}\frac{\varphi(u)}{u^2}\bigg(\sum_{\substack{n\in S(X)\\n\leq q^{1/10}/u\\(n,q)=1}}\frac{\alpha_{-1}(un)}{n}\bigg)^2.
\end{eqnarray*}
As before we can extend the sums to over all of $u,n\in S(X)$ with the gain of at most $O(q^{-1/50})$. Hence
\begin{equation*}
I_{1}^{M}=(1+o(1))\frac{\varphi(q)}{2q}\log q\sum_{\substack{u\in S(X)\\(u,q)=1}}\frac{\varphi(u)}{u^2}\bigg(\sum_{\substack{n\in S(X)\\(n,q)=1}}\frac{\alpha_{-1}(un)}{n}\bigg)^2+O(q^{-1/50}).
\end{equation*}
Since the functions $\varphi(n)$ and $\alpha_{-1}(n)$ are multiplicative, the entire sum is
\begin{equation*}
\prod_{\substack{p\leq X\\(p,q)=1}}\bigg(\sum_{h}\sum_{i}\sum_{j}\frac{\varphi(p^h)\alpha_{-1}(p^{h+i})\alpha_{-1}(p^{h+j})}{p^{2h+i+j}}\bigg).
\end{equation*}
We recall from \eqref{11} that for $p\leq\sqrt{X}$, $\alpha_{-1}(p)=-1$, $\alpha_{-1}(p^j)=0$ for every $j\geq2$, and for $\sqrt{X}<p\leq X$, $\alpha_{-1}(p)=1/2$, $\alpha_{-1}(p^j)=0$ for every $j\geq4$. Hence by Mertens' theorem, the above expression is
\begin{equation*}
\prod_{\substack{p\leq\sqrt{X}\\(p,q)=1}}\bigg(1-\frac{1}{p}\bigg)\prod_{\substack{\sqrt{X}<p\leq X\\(p,q)=1}}\bigg(1-\frac{1}{p}+O\bigg(\frac{1}{p^2}\bigg)\bigg)\sim\frac{1}{e^{\gamma}\log X}\prod_{\substack{p\leq X\\p|q}}\bigg(1-\frac{1}{p}\bigg)^{-1}.
\end{equation*}
So
\begin{equation*}
I_{1}^{M}=(1+o(1))\frac{1}{2}\prod_{\substack{p>X\\p|q}}\bigg(1-\frac{1}{p}\bigg)\frac{\log q}{e^{\gamma}\log X}+O\bigg(\prod_{p\leq X}\bigg(1+\frac{1}{p}\bigg)\bigg)+O(q^{-1/50}).
\end{equation*}
Thus
\begin{equation}\label{15}
I_1\sim\frac{1}{2}\prod_{\substack{p>X\\p|q}}\bigg(1-\frac{1}{p}\bigg)\frac{\log q}{e^{\gamma}\log X}.
\end{equation}

To estimate $I_{2}$, we note from Lemma 4 and the bounds for $W_{\mathfrak{a}}(x)$ that
\begin{equation*}
I_{2}\ll\frac{1}{\varphi^{*}(q)}\sum_{\substack{m,n\leq q^{1/10}\\(mn,q)=1}}\frac{|\alpha_{-1}(m)\alpha_{-1}(n)|}{\sqrt{mn}}\sum_{h|q}\varphi(h)\mu(q/h)^{2}\sum_{\substack{ab\leq Z\\am\equiv\pm bn(\textrm{mod}\ h)\\am\ne bn\\(ab,h)=1}}\frac{1}{\sqrt{ab}}.
\end{equation*}
Denote the innermost sum by $E_{1}(h)$. We divide the terms $ab\leq Z$ into dyadic blocks. Consider the block $Z_{1}\leq ab<2Z_{1}$. By Lemma 6, the sum over this block is
\begin{equation*}
\ll\frac{\sqrt{Z_{1}}mn}{h}\log (Z_{1}mn).
\end{equation*}
Summing over all dyadic blocks we have
\begin{equation*}
E_{1}(h)\ll\frac{\sqrt{Z}mn}{h}\log q.
\end{equation*}
So
\begin{equation*}
\frac{1}{\varphi^{*}(q)}\sum_{h|q}\varphi(h)\mu(q/h)^{2}E_{1}(h)\ll\frac{2^{\omega(q)/2}\sqrt{q}\log q}{\varphi^{*}(q)}mn\ll q^{-1/3}mn.
\end{equation*}
Thus
\begin{eqnarray}\label{16}
I_{2}&\ll&q^{-1/3}\sum_{m,n\leq q^{1/10}}|\alpha_{-1}(m)\alpha_{-1}(n)|\sqrt{mn}\nonumber\\
&\ll&q^{-1/3}\bigg(\sum_{m\leq q^{1/10}}m^{1/2+\epsilon}\bigg)^2
\ll q^{-1/30+\epsilon}.
\end{eqnarray}
This and \eqref{15} prove Proposition 1.
\end{proof}

\begin{prop}
Let $\delta>0$. Suppose we have $X,q\rightarrow\infty$ with $X\ll(\log q)^{2-\delta}$. Then
\begin{equation*}
J=\frac{1}{\varphi^{*}(q)}\sum_{\chi\ (\emph{mod}\ q)}{\!\!\!\!\!\!\!}^{\textstyle{*}}\ C(\chi)\bigg|\sum_{\substack{n\in S(X)\\n\leq q^{1/10}}}\frac{\alpha_{-1}(n)\chi(n)}{\sqrt{n}}\bigg|^2=o\bigg(\frac{\log q}{\log X}\bigg).
\end{equation*}
\end{prop}
\begin{proof}
We have
\begin{equation*}
J=\frac{1}{\varphi^{*}(q)}\sum_{\substack{ab>Z\\mn\in S(X)\\m,n\leq q^{1/10}}}\frac{\alpha_{-1}(m)\alpha_{-1}(n)}{\sqrt{abmn}}\sum_{\mathfrak{a}\in\{0,1\}}W_{\mathfrak{a}}\bigg(\frac{\pi ab}{q}\bigg)\sum_{\substack{\chi\ (\textrm{mod}\ q)\\\chi(-1)=(-1)^{\mathfrak{a}}}}{\!\!\!\!\!\!\!\!\!\!}^{\textstyle{*}}\ \chi(am)\overline{\chi}(bn)
\end{equation*}
We proceed as in Proposition 1. Let us write the last expression as $J_{1}+J_{2}$, where $J_{1}$ consists of the terms $am=bn$ and $J_{2}$ is the sum of the remaining terms. We first estimate $J_{1}$. In the same way as we dealt with $I_{1}$, we write $m=ur$, $n=us$, $a=vs$, and $b=vr$, where $(r,s)=1$. Then
\begin{eqnarray*}
J_{1}&=&\frac{1}{2}\sum_{\mathfrak{a}\in\{0,1\}}\sum_{\substack{v^{2}rs>Z\\urs\in S(X)\\ur,us\leq q^{1/10}\\(uvrs,q)=1\\(r,s)=1}}\frac{\alpha_{-1}(ur)\alpha_{-1}(us)}{uvrs}W_{\mathfrak{a}}\bigg(\frac{\pi v^{2}rs}{q}\bigg)\nonumber\\
&=&\frac{1}{2}\sum_{\mathfrak{a}\in\{0,1\}}\sum_{\substack{urs\in S(X)\\ur,us\leq q^{1/10}\\(urs,q)=1\\(r,s)=1}}\frac{\alpha_{-1}(ur)\alpha_{-1}(us)}{urs}\sum_{\substack{v>\sqrt{\frac{Z}{rs}} \\(v,q)=1}}\frac{1}{v}W_{\mathfrak{a}}\bigg(\frac{\pi v^{2}rs}{q}\bigg).
\end{eqnarray*}
The sum over $v$ is
\begin{eqnarray*}
&&\sum_{\substack{\sqrt{q/rs}>v>\sqrt{Z/rs}\\(v,q)=1}}\frac{1}{v}W_{\mathfrak{a}}\bigg(\frac{\pi v^{2}rs}{q}\bigg)+O\bigg(\sum_{v\geq\sqrt{q/rs}}\frac{1}{v}\bigg(\frac{v^{2}rs}{q}\bigg)^{-2}\bigg)\nonumber\\
&=&\sum_{\substack{\sqrt{q/rs}>v>\sqrt{Z/rs}\\(v,q)=1}}\frac{1}{v}+O\bigg(\sum_{v\leq\sqrt{q/rs}}\frac{1}{v}\frac{\sqrt{v}(rs)^{1/4}}{q^{1/4}}\bigg)+O(1)\nonumber\\
&=&\sum_{\substack{\sqrt{q/rs}>v>\sqrt{Z/rs}\\(v,q)=1}}\frac{1}{v}+O(1).
\end{eqnarray*}
By Lemma 5, this is
\begin{equation*}
\frac{\varphi(q)}{q}\big(\log\sqrt{q/rs}-\log\sqrt{Z/rs}+O(1+\log\omega(q))\big)+O(1)\ll\frac{\varphi(q)}{q}\omega(q),
\end{equation*}

Also as in Proposition 1, we have
\begin{eqnarray*}
\sum_{\substack{urs\in S(X)\\ur,us\leq q^{1/10}\\(urs,q)=1\\(r,s)=1}}\frac{|\alpha_{-1}(ur)\alpha_{-1}(us)|}{urs}\ll\log X.
\end{eqnarray*}
Hence
\begin{equation}\label{17}
J_{1}\ll \frac{\varphi(q)}{q}\omega(q)\log X=o\bigg(\frac{\log q}{\log X}\bigg).
\end{equation}

We now turn to $J_{2}$. We have
\begin{eqnarray*}
J_{2}\ll\frac{1}{\varphi^{*}(q)}\sum_{\substack{mn\in S(X)\\m,n\leq q^{1/10}\\(mn,q)=1}}\frac{|\alpha_{-1}(m)\alpha_{-1}(n)|}{\sqrt{mn}}\sum_{h|q}\varphi(h)\mu(q/h)^{2}\sum_{\substack{ab>Z\\am\equiv\pm bn(\textrm{mod}\ h)\\am\ne bn\\(ab,h)=1}}\frac{1}{\sqrt{ab}}\bigg(1+\frac{ab}{q}\bigg)^{-2}.
\end{eqnarray*}
We divide the innermost sum into dyadic blocks $Z_{1}<ab\leq2Z_{1}$, where $Z_{1}>Z$. We have
\begin{eqnarray*}
\sum_{\substack{Z_{1}<ab\leq2Z_{1}\\am\equiv\pm bn(\textrm{mod}\ h)\\am\ne bn\\(ab,h)=1}}\frac{1}{\sqrt{ab}}\bigg(1+\frac{ab}{q}\bigg)^{-2}&\ll&\frac{1}{\sqrt{Z_{1}}}\bigg(1+\frac{Z_{1}}{q}\bigg)^{-2}\sum_{\substack{Z_{1}<ab\leq2Z_{1}\\am\equiv\pm bn(\textrm{mod}\ h)\\am\ne bn\\(ab,h)=1}}1\nonumber\\
&\ll&\bigg(1+\frac{Z_{1}}{q}\bigg)^{-2}\frac{\sqrt{Z_{1}}mn}{h}\log (Z_{1}mn).\nonumber\\
&&
\end{eqnarray*}
Summing over all such blocks, we have that the innermost sum is $\ll\sqrt{q}mn(\log q)/h$. So, as for $I_{2}$,
\begin{equation}\label{18}
J_{2}\ll\frac{2^{\omega(q)}\sqrt{q}}{\varphi^{*}(q)\log q}\sum_{m,n\leq q^{1/10}}|\alpha_{-1}(m)\alpha_{-1}(n)|\sqrt{mn}\ll q^{-1/30+\epsilon}.
\end{equation}
The proof of Proposition 2 is complete.
\end{proof}

The first part of Theorem 3 follows from Proposition 1, Proposition 2 and \eqref{12} with $\theta = 1/10$.

The second statement then follows by Lemma 2.

\section{Proof of Theorem 4}

We recall from Lemma 3 that
\begin{equation}\label{19}
P_{X}(\chi)^{-2}=\sum_{n\in S(X)}\frac{\alpha_{-2}(n)\chi(n)}{\sqrt{n}}\bigg(1+O_{k}\bigg(\frac{1}{\log X}\bigg)\bigg),
\end{equation}
where $\alpha_{-2}(n)$ is defined by
\begin{equation*}
\sum_{n=1}^{\infty}\frac{\alpha_{-2}(n)\chi(n)}{n^s}=\prod_{p\leq X}\bigg( 1-\frac{\chi(p)}{p^s}\bigg) ^{2}\prod_{\sqrt{X}<p\leq X}\bigg( 1+\frac{\chi(p)^2}{2p^{2s}}\bigg) ^{2}.
\end{equation*}
We note that $\alpha_{-2}(p^l)=0$ for $l\geq7$, and for $l\geq3$ if $p\leq\sqrt{X}$. 

We first need a lemma relating to \eqref{19}.

\begin{lemm}
Let $\beta_{-2}(n)$ be a multiplicative function defined by 
\begin{equation*}
\beta_{-2}(p)=\alpha_{-2}(p),\quad \beta_{-2}(p^2)=\alpha_{-2}(p^2),\quad \emph{and}\quad \beta_{-2}(p^l)=0\ \emph{for}\ l\geq3.
\end{equation*}
Then
\begin{equation*}
\sum_{n\in S(X)}\frac{\alpha_{-2}(n)\chi(n)}{\sqrt{n}}=\bigg(\sum_{n\in S(X)}\frac{\beta_{-2}(n)\chi(n)}{\sqrt{n}}\bigg)\bigg(1+O\bigg(\frac{1}{X^{1/4}\log X}\bigg)\bigg).
\end{equation*}
\end{lemm}
\begin{proof}
Define by $S'(X)$ the subset of $S(X)$ consisting of cube-free integers. We have
\begin{eqnarray*}
\sum_{n\in S(X)}\frac{\alpha_{-2}(n)\chi(n)}{\sqrt{n}}&=&\sum_{n\in S'(X)}\frac{\beta_{-2}(n)\chi(n)}{\sqrt{n}}\prod_{\substack{\sqrt{X}<p\leq X\\p\nmid n}}\bigg(1+\sum_{l\geq3}\frac{\alpha_{-2}(p^l)\chi(p^l)}{p^{l/2}}\bigg)\nonumber\\
&=&\sum_{n\in S'(X)}\frac{\beta_{-2}(n)\chi(n)}{\sqrt{n}}\exp\bigg(O\bigg(\sum_{\sqrt{X}<p\leq X}\frac{1}{p^{3/2}}\bigg)\bigg)\nonumber\\
&=&\sum_{n\in S'(X)}\frac{\beta_{-2}(n)\chi(n)}{\sqrt{n}}\bigg(1+O\bigg(\frac{1}{X^{1/4}\log X}\bigg)\bigg).
\end{eqnarray*}
Since $\beta_{-2}(n)=0$ for $n\in S(X)-S'(X)$, the lemma follows.
\end{proof}

\textbf{Remark}.\quad Lemma 7 implies that we may assume $\alpha_{-2}(n)$ is supported on cube-free integers.\\

Similar to \eqref{12}, for $X\ll(\log\log q)^{2-\delta}$, we can truncate the series for $P_{X}^{*}(\chi)^{-2}$ at $(\log q)^{1/4}$. We have
\begin{equation*}
\sum_{n\in S(X)}\frac{\alpha_{-2}(n)\chi(n)}{\sqrt{n}}=\sum_{\substack{n\in S(X)\\n\leq (\log q)^{1/4}}}\frac{\alpha_{-2}(n)\chi(n)}{\sqrt{n}}+O\bigg(\sum_{\substack{n\in S(X)\\n>(\log q)^{1/4}}}\frac{d_3(n)}{\sqrt{n}}\bigg).
\end{equation*}
The $O$-term is
\begin{eqnarray*}
&\ll&\sum_{n\in S(X)}\bigg(\frac{n}{(\log q)^{1/4}}\bigg)^{\delta/4}\frac{d_3(n)}{\sqrt{n}}=(\log q)^{-\delta/16}\prod_{p\leq X}\bigg(1-\frac{1}{p^{1/2-\delta/4}}\bigg)^{-3}\\
&\ll&(\log q)^{-\delta/16}\exp\bigg(O\bigg(\frac{X^{1/2+\delta/4}}{\log X}\bigg)\bigg)\\
&\ll&(\log q)^{-\delta/16}\exp\bigg(O\bigg(\frac{\log\log q}{\log\log\log q}\bigg)\bigg)\ll(\log q)^{-\delta/17},
\end{eqnarray*}
since $X^{1/2+\delta/4}\ll(\log\log q)^{(2-\delta)(1/2+\delta/4)}\ll\log\log q$. Thus
\begin{equation}\label{20}
P_{X}^{*}(\chi)^{-2}=\sum_{\substack{n\in S(X)\\n\leq (\log q)^{1/4}}}\frac{\alpha_{-2}(n)\chi(n)}{\sqrt{n}}+O((\log q)^{-\delta/17}).
\end{equation}

\begin{lemm}
For $x\geq\sqrt{q}$, $h\leq q$ and $l=0$ or $1$, we have
\begin{equation*}
\sum_{\substack{n\leq x\\(n,q)=1}}\frac{2^{\omega(n)-\omega((n,h))}}{n}\bigg(\log\frac{x}{n}\bigg)^l\ll(\log x)^{2+l},
\end{equation*}
and
\begin{equation*}
\sum_{\substack{n\leq x\\(n,q)=1}}\frac{2^{\omega(n)-\omega((n,h))}}{n}\bigg(\log\frac{x}{n}\bigg)^2\sim\frac{(\log x)^{4}}{12\zeta(2)}\prod_{p|hq}\bigg(\frac{1-1/p}{1+1/p}\bigg)\prod_{\substack{p|h\\p\nmid q}}\bigg(\frac{1}{1-1/p}\bigg).
\end{equation*}
\end{lemm}
\begin{proof}
For $\Re s>1$, we define
\begin{equation*}
F(s)=\sum_{\substack{n=1\\(n,q)=1}}^{\infty}\frac{2^{\omega(n)-\omega((n,h))}}{n^s}=\frac{\zeta(s)^2}{\zeta(2s)}\prod_{p|hq}\frac{1-p^{-s}}{1+p^{-s}}\prod_{\substack{p|h\\p\nmid q}}\frac{1}{1-p^{-s}}.
\end{equation*}
We have
\begin{equation}\label{28}
\sum_{\substack{n\leq x\\(n,q)=1}}\frac{2^{\omega(n)-\omega((n,h))}}{n}\ll\sum_{\substack{n=1\\(n,q)=1}}^{\infty}\frac{2^{\omega(n)-\omega((n,h))}}{n^{1+1/\log x}}\ll F(1+1/\log x)\ll (\log x)^2.
\end{equation}
We obtain the result for $l=0$. Furthermore, we note that
\begin{equation*}
\frac{1}{2}\sum_{\substack{n\leq x\\(n,q)=1}}\frac{2^{\omega(n)-\omega((n,h))}}{n}\bigg(\log\frac{x}{n}\bigg)^2=\frac{1}{2\pi i}\int_{(c)}F(s+1)\frac{x^{s}}{s^3}ds,
\end{equation*}
for some $c>0$. By moving the line of integration to $-1/2+\epsilon$ we obtain
\begin{equation*}
\frac{1}{2}\sum_{\substack{n\leq x\\(n,q)=1}}\frac{2^{\omega(n)-\omega((n,h))}}{n}\bigg(\log\frac{x}{n}\bigg)^2=\textrm{Res}_{s=0}F(s+1)\frac{x^s}{s^3}+O(q^{\epsilon}x^{-1/2+\epsilon}).
\end{equation*}
Also, it is easy to check that
\begin{equation*}
\textrm{Res}_{s=0}F(s+1)\frac{x^s}{s^3}=\frac{(\log x)^{4}}{24\zeta(2)}\prod_{p|hq}\bigg(\frac{1-1/p}{1+1/p}\bigg)\prod_{\substack{p|h\\p\nmid q}}\bigg(\frac{1}{1-1/p}\bigg)\bigg(1+O\bigg(\frac{1+\log\omega(hq)}{\log q}\bigg)\bigg).
\end{equation*}
The last statement follows. The case $l=1$ follows by combining this with \eqref{28} and Cauchy's inequality.
\end{proof}

We mention a result of Shiu [\textbf{\ref{Sh}}].

\begin{lemm}
If $(r,h)=1$ and $h^{1+\delta}\leq x$ for some $\delta>0$ then
\begin{equation*}
\sum_{\substack{n\leq x\\n\equiv r(\emph{mod}\ h)}}d(n)\ll\frac{\varphi(h)}{h^2}x\log x.
\end{equation*}
\end{lemm}

We next mimic a result of Heath-Brown [\textbf{\ref{H-B2}}] to give the following result.

\begin{lemm}
For $l$, $h$, $m$ positive integers which satisfy $(lh)^{4/5}\ll x$, we have
\begin{equation*}
\sum_{\substack{n\leq x\\(n,h)=1}}d(n)d(lh+mn)\ll mx(\log mx)^{2}\sum_{d|l}d^{-1},
\end{equation*}
and if we also have $mx<lh$ then
\begin{equation*}
\sum_{\substack{n\leq x\\(n,h)=1}}d(n)d(lh-mn)\ll mx(\log mx)^{2}\sum_{d|l}d^{-1}.
\end{equation*}
\end{lemm}
\begin{proof}
We prove the first part. The second part can be done similarly. 

We use an estimate given by Heath-Brown in [\textbf{\ref{H-B1}}] which asserts that given $kX^{2/5}\ll x$, we have
\begin{equation*}
\sum_{\substack{X-x<n\leq X\\n\equiv l(\textrm{mod}\ k)}}d(n)\ll C(k,l)x\log x,
\end{equation*}
where
\begin{equation*}
C(k,l)=\sum_{\substack{d|(k,l)\\\delta|kd^{-1}}}d\delta k^{-2}.
\end{equation*}
For $n\leq x$ and $(n,h)=1$ we have
\begin{equation*}
d(n)\ll\sum_{\substack{k\leq\sqrt{x}\\k|n\\(k,h)=1}}1.
\end{equation*}
So
\begin{eqnarray*}
\sum_{\substack{n\leq x\\(n,h)=1}}d(n)d(lh+mn)&\ll&\sum_{\substack{k\leq\sqrt{x}\\(k,h)=1}}\sum_{\substack{lh<u\leq lh+mx\\u\equiv lh(\textrm{mod}\ k)}}d(u)\nonumber\\
&\ll&\sum_{k\leq\sqrt{x}}C(k,l)mx\log mx\ll mx\log mx\sum_{k\leq\sqrt{x}}\sum_{\substack{d|(k,l)\\\delta|kd^{-1}}}d\delta k^{-2}\nonumber\\
&\ll&mx\log mx\sum_{\substack{d\delta\leq\sqrt{x}\\d|l}}(d\delta)^{-1}\ll mx(\log mx)^{2}\sum_{d|l}d^{-1},
\end{eqnarray*}
since $k(lh+mx)^{2/5}\ll x^{1/2}((lh)^{2/5}+(mx)^{2/5})\ll mx$. This completes the proof.
\end{proof}

\begin{lemm}
Let $m$, $n$, $h$ be three positive integers, $(mn,h)=1$ and $Z_{1}, Z_{2}\geq 2$. If $Z_{1}Z_{2}>h^{19/10}$ then we have
\begin{equation*}
E_{2}:=\sum_{\substack{Z_{1}\leq ab<2Z_{1}\\Z_{2}\leq cd<2Z_{2}\\acm\equiv\pm bdn(\emph{mod}\ h)\\acm\ne bdn\\(abcd,h)=1}}1\ll\frac{Z_{1}Z_{2}mn}{h}(\log Z_{1}Z_{2}mn)^3,
\end{equation*}
and if $Z_{1}Z_{2}\leq h^{19/10}$ then
\begin{equation*}
E_{2}\ll\frac{(Z_{1}Z_{2}mn)^{1+\epsilon}}{h}.
\end{equation*}
\end{lemm}
\begin{proof}
It is clear that we only need to consider the case $acm>bdn$. Let $acm=lh\pm bdn$, where $l\geq1$, and $u=bd$. We have 
\begin{equation*}
l\leq\frac{acm+bdn}{h}\leq\frac{4Z_{1}Z_{2}(m+n)}{h},\quad \textrm{and}\quad u\leq\frac{4Z_{1}Z_{2}}{ac}\leq\frac{8Z_{1}Z_{2}m}{lh}.
\end{equation*}
So
\begin{eqnarray}\label{21}
E_{2}&\ll&\sum_{1\leq l\leq\frac{4Z_{1}Z_{2}(m+n)}{h}}\sum_{\substack{u\leq\frac{8Z_{1}Z_{2}m}{lh}\\lh\pm nu\equiv0(\textrm{mod}\ m)\\(u,h)=1}}d(u)d\bigg(\frac{lh\pm nu}{m}\bigg)\nonumber\\
&\ll&\sum_{1\leq l\leq\frac{4Z_{1}Z_{2}(m+n)}{h}}\sum_{\substack{u\leq\frac{8Z_{1}Z_{2}m}{lh}\\(u,h)=1}}d(u)d(lh\pm nu).
\end{eqnarray}

We consider first the case $Z_{1}Z_{2}>h^{19/10}$. If we also have $lh\leq(Z_{1}Z_{2})^{5/9}$ then by Lemma 10,
\begin{eqnarray*}
E_{2}&\ll&\sum_{1\leq l\leq\frac{(Z_{1}Z_{2})^{5/9}}{h}}\frac{Z_{1}Z_{2}mn}{lh}(\log Z_{1}Z_{2}mn)^{2}\sum_{d|l}\frac{1}{d}\nonumber\\
&\ll&\frac{Z_{1}Z_{2}mn}{h}(\log Z_{1}Z_{2}mn)^{3},
\end{eqnarray*}
as $8Z_{1}Z_{2}m/lh>Z_{1}Z_{2}/lh\gg(lh)^{4/5}$. For $lh>(Z_{1}Z_{2})^{5/9}$, let $v=acm$. Then
\begin{equation*}
v>lh/2>(Z_{1}Z_{2})^{5/9}/2,\quad \textrm{and}\quad v\leq4Z_{1}Z_{2}m/bd=4Z_{1}Z_{2}m/u.
\end{equation*}
We note that $d(ac)\leq d(v)$. Hence the expression in \eqref{21} is
\begin{equation*}
\sum_{\substack{u\leq\frac{8Z_{1}Z_{2}m}{h}\\(u,h)=1}}d(u)\sum_{\substack{\frac{(Z_{1}Z_{2})^{5/9}}{2}<v\leq\frac{4Z_{1}Z_{2}m}{u}\\v\equiv\pm un(\textrm{mod}\ h)}}d(v),
\end{equation*}
which is, by Lemma 9,
\begin{eqnarray*}
&\ll&\sum_{u\leq\frac{8Z_{1}Z_{2}m}{h}}d(u)\frac{\varphi(h)}{h^2}\frac{Z_{1}Z_{2}m}{u}\log Z_{1}Z_{2}m\ll\frac{Z_{1}Z_{2}mn}{h}(\log Z_{1}Z_{2}mn)^{3},
\end{eqnarray*}
since $(Z_{1}Z_{2})^{5/9}/2\gg h^{95/90}$. This proves the first statement.

For $Z_{1}Z_{2}\leq h^{19/10}$, since $d(u)d(lh\pm nu)\ll(Z_{1}Z_{2}mn)^{\epsilon}$, from \eqref{21} we have
\begin{eqnarray*}
E_{2}\ll\sum_{1\leq l\leq\frac{4Z_{1}Z_{2}(m+n)}{h}}\frac{Z_{1}Z_{2}m}{lh}(Z_{1}Z_{2}mn)^{\epsilon}\ll\frac{(Z_{1}Z_{2}mn)^{1+\epsilon}}{h}.
\end{eqnarray*}
The proof of the lemma is complete.
\end{proof}

We are now ready to prove two propositions.

\begin{prop}
Let $\delta>0$. Suppose we have $X,q\rightarrow\infty$ with $X\ll(\log\log q)^{2-\delta}$. Then
\begin{equation*}
I=\frac{1}{\varphi^{*}(q)}\sum_{\chi\ (\emph{mod}\ q)}{\!\!\!\!\!\!\!}^{\textstyle{*}}\ \ \bigg|B(\chi)\sum_{\substack{n\in S(X)\\n\leq (\log q)^{1/4}}}\frac{\alpha_{-2}(n)\chi(n)}{\sqrt{n}}\bigg|^2\sim\frac{1}{48}\prod_{\substack{p>X\\p|q}}\frac{(1-1/p)^{3}}{1+1/p}\bigg(\frac{\log q}{e^{\gamma}\log X}\bigg)^{4}.
\end{equation*}
\end{prop}
\begin{proof}
We have
\begin{eqnarray*}
I&=&\frac{1}{\varphi^{*}(q)}\sum_{\substack{ab,cd\leq Z\\mn\in S(X)\\m,n\leq (\log q)^{1/4}}}\frac{\alpha_{-2}(m)\alpha_{-2}(n)}{\sqrt{abcdmn}}\sum_{\mathfrak{a}\in\{0,1\}}W_{\mathfrak{a}}\bigg(\frac{\pi ab}{q}\bigg)W_{\mathfrak{a}}\bigg(\frac{\pi cd}{q}\bigg)\sum_{\substack{\chi\ (\textrm{mod}\ q)\\\chi(-1)=(-1)^{\mathfrak{a}}}}{\!\!\!\!\!\!\!\!\!}^{\textstyle{*}}\ \chi(acm)\overline{\chi}(bdn)\nonumber\\
&=&I_{1}+I_{2},
\end{eqnarray*}
where $I_{1}$ and $I_{2}$ are, respectively, the diagonal and the off-diagonal. We first consider $I_1$:
\begin{equation*}
I_{1}=\frac{1}{2}\sum_{\mathfrak{a}\in\{0,1\}}\sum_{\substack{ab,cd\leq Z\\mn\in S(X)\\m,n\leq (\log q)^{1/4}\\acm=bdn\\(abcdmn,q)=1}}\frac{\alpha_{-2}(m)\alpha_{-2}(n)}{\sqrt{abcdmn}}W_{\mathfrak{a}}\bigg(\frac{\pi ab}{q}\bigg)W_{\mathfrak{a}}\bigg(\frac{\pi cd}{q}\bigg).
\end{equation*}
Since $acm=bdn$, we can write $m=ugh$, $n=uij$, $a=vjk$, $b=vgl$, $c=wil$ and $d=whk$, where $(gh,ij)=(k,l)=(k,gi)=(l,hj)=1$. Also let $f=kl$. We note that, given $f$, there are $2^{\omega(f)-\omega((f,ghij))}$ ways to express $f$ as $kl$ such that $(k,l)=(k,gi)=(l,hj)=1$. Hence
\begin{eqnarray*}
I_{1}&=&\frac{1}{2}\sum_{\mathfrak{a}\in\{0,1\}}\sum_{\substack{v^{2}gjf,w^{2}hif\leq Z\\ughij\in S(X)\\ugh,uij\leq (\log q)^{1/4}\\(uvwghij,q)=1\\(gh,ij)=1\\(f,q(gi,hj))=1}}\frac{2^{\omega(f)-\omega((f,ghij))}\alpha_{-2}(ugh)\alpha_{-2}(uij)}{uvwghijf}W_{\mathfrak{a}}\bigg(\frac{\pi v^{2}gjf}{q}\bigg)W_{\mathfrak{a}}\bigg(\frac{\pi w^{2}hif}{q}\bigg)\nonumber\\
&=&\frac{1}{2}\sum_{\mathfrak{a}\in\{0,1\}}\sum_{\substack{ughij\in S(X)\\ugh,uij\leq (\log q)^{1/4}\\(ughij,q)=1\\(gh,ij)=1}}\frac{\alpha_{-2}(ugh)\alpha_{-2}(uij)}{ughij}\sum_{\substack{f\leq\min\{Z/gj,Z/hi\} \\(f,q(gi,hj))=1}}\frac{2^{\omega(f)-\omega((f,ghij))}}{f}\nonumber\\
&&\quad\bigg[\sum_{\substack{v\leq\sqrt{Z/gjf}\\(v,q)=1}}\frac{1}{v}W_{\mathfrak{a}}\bigg(\frac{\pi v^{2}gjf}{q}\bigg)\bigg]\bigg[\sum_{\substack{w\leq\sqrt{Z/hif}\\(w,q)=1}}\frac{1}{w}W_{\mathfrak{a}}\bigg(\frac{\pi w^{2}hif}{q}\bigg)\bigg].
\end{eqnarray*}
The expression in the first square bracket is
\begin{equation*}
\sum_{\substack{v\leq\sqrt{Z/gjf}\\(v,q)=1}}\frac{1}{v}\bigg(1+O\bigg(\frac{\sqrt{v}(gjf)^{1/4}}{q^{1/4}}\bigg)\bigg)=\sum_{\substack{v\leq\sqrt{Z/gjf}\\(v,q)=1}}\frac{1}{v}+O(2^{-\omega(q)/4}),
\end{equation*}
which is, by Lemma 5, equal to
\begin{equation*}
\frac{\varphi(q)}{2q}\log\frac{Z}{gjf}+O(1+\log\omega(q))=\frac{\varphi(q)}{2q}\log\frac{Z}{f}+O(\log\log q),
\end{equation*}
if $gjf\leq Z_{0}=Z/9^{\omega(q)}$, and if $gjf>Z_{0}$,
\begin{equation*}
\ll\sum_{\substack{v\leq3^{\omega(q)}\\(v,q)=1}}\frac{1}{v}+O(2^{-\omega(q)/4})\ll\omega(q).
\end{equation*}
A similar calculation holds for the second square bracket, so, for $m_{0}=\min\{Z_{0}/gj,Z_{0}/hi\}$ and $M_{0}=\max\{Z_{0}/gj,Z_{0}/hi\}$, the sum over $f$ is $R_{1}+R_{2}+R_{3}$, where, by Lemma 8,
\begin{equation*}
R_{3}\ll\omega(q)^{2}\sum_{\substack{M_{0}\leq f\leq\min\{Z/gj,Z/hi\}\\(f,q(gi,hj))=1}}\frac{2^{\omega(f)-\omega((f,ghij))}}{f}\ll\omega(q)^{2}(\log q)^2,
\end{equation*}
\begin{eqnarray*}
R_{2}&\ll&\omega(q)\sum_{\substack{m_{0}\leq f\leq M_{0}\\(f,q(gi,hj))=1}}\frac{2^{\omega(f)-\omega((f,ghij))}}{f}\bigg(\log\frac{Z}{f}+O(\log\log q)\bigg)\nonumber\\
&\ll&\omega(q)(\log q)^{3},
\end{eqnarray*}
and
\begin{eqnarray*}
R_{1}&=&\frac{\varphi(q)^2}{4q^2}\sum_{\substack{f\leq m_{0}\\(f,q(gi,hj))=1}}\frac{2^{\omega(f)-\omega((f,ghij))}}{f}\bigg(\log\frac{Z}{f}+O(\log\log q)\bigg)^2\nonumber\\
&\sim&\frac{\varphi(q)^2}{4q^2}\frac{(\log q)^4}{12\zeta(2)}\prod_{p|qghij}\frac{1-1/p}{1+1/p}\prod_{\substack{p|ghij\\p\nmid (gi,hj)}}\frac{1}{1-1/p}\nonumber\\
&\sim&\frac{1}{8\pi^{2}}\prod_{p|q}\frac{(1-1/p)^3}{1+1/p}(\log q)^4\prod_{p|ghij}\frac{1-1/p}{1+1/p}\prod_{\substack{p|ghij\\p\nmid (g,h)(i,j)}}\frac{1}{1-1/p}.
\end{eqnarray*}
Thus
\begin{eqnarray}\label{22}
I_{1}&\sim&\frac{1}{8\pi^{2}}\prod_{p|q}\frac{(1-1/p)^3}{1+1/p}(\log q)^4\nonumber\\
&&\quad\times\sum_{\substack{ughij\in S(X)\\ugh,uij\leq (\log q)^{1/4}\\(ughij,q)=1\\(gh,ij)=1}}\frac{\alpha_{-2}(ugh)\alpha_{-2}(uij)}{ughij}\prod_{p|ghij}\frac{1-1/p}{1+1/p}\prod_{\substack{p|ghij\\p\nmid (g,h)(i,j)}}\frac{1}{1-1/p}.
\end{eqnarray}
We need to estimate the last factor, which is
\begin{equation*}
S=\sum_{\substack{u\in S(X)\\u\leq (\log q)^{1/4}\\(u,q)=1}}\frac{1}{u}\sum_{\substack{m\in S(X)\\m\leq (\log q)^{1/4}/u\\(m,q)=1}}\sum_{\substack{n\in S(X)\\n\leq (\log q)^{1/4}/u\\(n,mq)=1}}\frac{\alpha_{-2}(um)\alpha_{-2}(un)\delta(m)\delta(n)}{mn},
\end{equation*}
where
\begin{equation}
\delta(m)=\prod_{p^{r}||m}\frac{1-1/p}{1+1/p}\bigg(\frac{2}{1-1/p}+(r-1)\bigg)=\prod_{p^{r}||m}\bigg(1+r\frac{1-1/p}{1+1/p}\bigg).
\end{equation}
Let $P=\prod_{p\leq X}p$. Since $\alpha_{-2}(n)=0$ if $n$ is not a cube-free integer, we can restrict the summation over $u$ to summation over $u=u_{1}u_{2}^{2}$, where $u_{1}|P$, and $u_{2}|(P/u_{1})$. The summation over $m$ and $n$ can also be restricted to $(m,u_{2})=(n,u_{2})=1$, since otherwise $\alpha_{-2}(um)\alpha_{-2}(un)=0$. So
\begin{eqnarray*}
&&S=\sum_{\substack{u_{1}\leq (\log q)^{1/4}\\u_{1}|P\\(u_{1},q)=1}}\frac{1}{u_{1}}\sum_{\substack{u_{2}\leq\sqrt{(\log q)^{1/4}/u_{1}}\\u_{2}|(P/u_{1})\\(u_{2},q)=1}}\frac{\alpha_{-2}(u_{2}^{2})^2}{u_{2}^{2}}\nonumber\\
&&\qquad\qquad\qquad\sum_{\substack{m\in S(X)\\m\leq (\log q)^{1/4}/u_{1}u_{2}^{2}\\(m,qu_{2})=1}}\sum_{\substack{n\in S(X)\\n\leq (\log q)^{1/4}/u_{1}u_{2}^{2}\\(n,mqu_{2})=1}}\frac{\alpha_{-2}(u_{1}m)\alpha_{-2}(u_{1}n)\delta(m)\delta(n)}{mn}.
\end{eqnarray*}
The same arguments imply that if $r=(u_{1},m)$ and $m=rm_{1}$ then we can assume that $(r,m_{1})=1$, and hence $(u_{1},m_{1})=1$. The summation over $m$ can be replaced by
\begin{equation*}
\sum_{r|u_{1}}\sum_{\substack{m_{1}\in S(X)\\m_{1}\leq (\log q)^{1/4}/u_{1}u_{2}^{2}r\\(m_{1},qu_{1}u_{2})=1}}.
\end{equation*}
Similarly, for $s=(u_{1},n)$ and $n=sn_{1}$ we can sum over $(u_{1},n_{1})=1$. The condition $(m,n)=1$ is equivalent to $(m_{1},n_{1})=(m_{1},s)=(r,n_{1})=(r,s)=1$. We have $(r,s)=1$ if and only if $s|(u_{1}/r)$. Also $(m_{1},s)=1$ and $(n_1,r)=1$ automatically come from $(m_{1}n_{1},u_{1})=1$. Thus
\begin{eqnarray*}
S&=&\sum_{\substack{u_{1}\leq (\log q)^{1/4}\\u_{1}|P\\(u_{1},q)=1}}\frac{1}{u_{1}}\sum_{\substack{u_{2}\leq\sqrt{(\log q)^{1/4}/u_{1}}\\u_{2}|(P/u_{1})\\(u_{2},q)=1}}\frac{\alpha_{-2}(u_{2}^{2})^2}{u_{2}^{2}}\sum_{r|u_{1}}\sum_{\substack{m_{1}\in S(X)\\m_{1}\leq (\log q)^{1/4}/u_{1}u_{2}^{2}r\\(m_{1},qu_{1}u_{2})=1}}\nonumber\\
&&\qquad\qquad\sum_{s|(u_{1}/r)}\sum_{\substack{n_{1}\in S(X)\\n_{1}\leq (\log q)^{1/4}/u_{1}u_{2}^{2}s\\(n_{1},m_{1}qu_{1}u_{2})=1}}\frac{\alpha_{-2}(u_{1}rm_{1})\alpha_{-2}(u_{1}sn_{1})\delta(rm_{1})\delta(sn_{1})}{rsm_{1}n_{1}}\nonumber\\
&=&\sum_{\substack{u_{1}\leq (\log q)^{1/4}\\u_{1}|P\\(u_{1},q)=1}}\frac{\alpha_{-2}(u_{1})^2}{u_{1}}\sum_{\substack{u_{2}\leq\sqrt{(\log q)^{1/4}/u_{1}}\\u_{2}|(P/u_{1})\\(u_{2},q)=1}}\frac{\alpha_{-2}(u_{2}^{2})^2}{u_{2}^{2}}\nonumber\\
&&\qquad\qquad\qquad\qquad\qquad\sum_{r|u_{1}}\frac{\alpha_{-2}(r^2)\delta(r)}{\alpha_{-2}(r)r}\sum_{s|(u_{1}/r)}\frac{\alpha_{-2}(s^2)\delta(s)}{\alpha_{-2}(s)s}\nonumber\\
&&\quad\sum_{\substack{m_{1}\in S(X)\\m_{1}\leq (\log q)^{1/4}/u_{1}u_{2}^{2}r\\(m_{1},qu_{1}u_{2})=1}}\frac{\alpha_{-2}(m_{1})\delta(m_{1})}{m_{1}}\sum_{\substack{n_{1}\in S(X)\\n_{1}\leq (\log q)^{1/4}/u_{1}u_{2}^{2}s\\(n_{1},m_{1}qu_{1}u_{2})=1}}\frac{\alpha_{-2}(n_{1})\delta(n_{1})}{n_{1}}.
\end{eqnarray*}
We can extend all the sums to all of $S(X)$ and $(u_1u_2m_1n_1,q)=1$ with the gain of at most `little $o$' of the main term. Also note that we can just consider $m_1$ and $n_1$ to be cube-free so we have
\begin{eqnarray*}
S&=&(1+o(1))\sum_{\substack{u_{1}|P\\(u_{1},q)=1}}\frac{\alpha_{-2}(u_{1})^2}{u_{1}}\sum_{\substack{u_{2}|(P/u_{1})\\(u_{2},q)=1}}\frac{\alpha_{-2}(u_{2}^{2})^2}{u_{2}^{2}}\sum_{r|u_{1}}\frac{\alpha_{-2}(r^2)\delta(r)}{\alpha_{-2}(r)r}\sum_{s|(u_{1}/r)}\frac{\alpha_{-2}(s^2)\delta(s)}{\alpha_{-2}(s)s}\nonumber\\
&&\qquad\sum_{\substack{m_{1}|(P/u_{1}u_{2})^2\\(m_{1},q)=1}}\frac{\alpha_{-2}(m_{1})\delta(m_{1})}{m_{1}}\sum_{\substack{n_{1}|(P/m_{1}u_{1}u_{2})^{2}\\(n_{1},q)=1}}\frac{\alpha_{-2}(n_{1})\delta(n_{1})}{n_{1}}.
\end{eqnarray*}
Define the following multiplicative functions
\begin{eqnarray*}
&&T_{1}(n):=\sum_{\substack{h|n\\(h,q)=1}}\frac{\alpha_{-2}(h)\delta(h)}{h},\ T_{2}(n):=\sum_{\substack{h|n\\(h,q)=1}}\frac{\alpha_{-2}(h)\delta(h)}{hT_{1}(h^2)},\\
&&T_{3}(n):=\sum_{\substack{h|n\\(h,q)=1}}\frac{\alpha_{-2}(h^2)\delta(h)}{\alpha_{-2}(h)h},\ T_{4}(n):=\sum_{\substack{h|n\\(h,q)=1}}\frac{\alpha_{-2}(h^2)\delta(h)}{\alpha_{-2}(h)hT_{3}(h)},\\
&&T_{5}(n):=\sum_{\substack{h|n\\(h,q)=1}}\frac{\alpha_{-2}(h^2)^{2}}{h^{2}T_{1}(h^2)T_{2}(h^2)},\ T_{6}(n):=\sum_{\substack{h|n\\(h,q)=1}}\frac{\alpha_{-2}(h)^{2}T_{3}(h)T_{4}(h)}{hT_{1}(h^2)T_{2}(h^2)T_{5}(h)}.
\end{eqnarray*}
The innermost sum of $S$ is
\begin{equation*}
T_{1}((P/m_{1}u_{1}u_{2})^2)=T_{1}(P^2)/T_{1}(m_{1}^{2})T_{1}(u_{1}^{2})T_{1}(u_{2}^{2}).
\end{equation*}
So the contribution of summing over $m_{1}$ and $n_{1}$ is
\begin{equation*}
\frac{T_{1}(P^2)T_{2}(P^2)}{T_{1}(u_{1}^{2})T_{1}(u_{2}^{2})T_{2}(u_{1}^{2})T_{2}(u_{2}^{2})}.
\end{equation*}
Thus
\begin{eqnarray*}
S&=&(1+o(1))T_{1}(P^2)T_{2}(P^2)\sum_{\substack{u_{1}|P\\(u_{1},q)=1}}\frac{\alpha_{-2}(u_{1})^2}{u_{1}T_{1}(u_{1}^{2})T_{2}(u_{1}^{2})}\sum_{\substack{u_{2}|(P/u_{1})\\(u_{2},q)=1}}\frac{\alpha_{-2}(u_{2}^{2})^2}{u_{2}^{2}T_{1}(u_{1}^{2})T_{2}(u_{1}^{2})}\nonumber\\
&&\qquad\qquad\sum_{r|u_{1}}\frac{\alpha_{-2}(r^2)\delta(r)}{\alpha_{-2}(r)r}\sum_{s|(u_{1}/r)}\frac{\alpha_{-2}(s^2)\delta(s)}{\alpha_{-2}(s)s}\nonumber\\
&=&(1+o(1))T_{1}(P^2)T_{2}(P^2)\sum_{\substack{u_{1}|P\\(u_{1},q)=1}}\frac{\alpha_{-2}(u_{1})^2T_{3}(u_{1})T_{4}(u_{1})}{u_{1}T_{1}(u_{1}^{2})T_{2}(u_{1}^{2})}\sum_{\substack{u_{2}|(P/u_{1})\\(u_{2},q)=1}}\frac{\alpha_{-2}(u_{2}^{2})^2}{u_{2}^{2}T_{1}(u_{1}^{2})T_{2}(u_{1}^{2})}\nonumber\\
&=&(1+o(1))T_{1}(P^2)T_{2}(P^2)T_{5}(P)T_{6}(P)\nonumber\\
&=&(1+o(1))\prod_{\substack{p|P\\p\nmid q}}\bigg(T_{1}(p^2)T_{2}(p^2)T_{5}(p)+\frac{\alpha_{-2}(p)^{2}T_{3}(p)T_{4}(p)}{p}\bigg).
\end{eqnarray*}
Since $\alpha_{-2}(p)=-2$,
\begin{eqnarray*}
T_{3}(p)T_{4}(p)&=&T_{3}(p)+\frac{\alpha_{-2}(p^2)\delta(p)}{\alpha_{-2}(p)p}\nonumber\\
&=&1+\frac{\alpha_{-2}(p^2)\delta(p)}{\alpha_{-2}(p)p}+\frac{\alpha_{-2}(p^2)\delta(p)}{\alpha_{-2}(p)p}\nonumber\\
&=&1-\frac{\alpha_{-2}(p^2)\delta(p)}{p},
\end{eqnarray*}
and
\begin{eqnarray*}
T_{1}(p^2)T_{2}(p^2)T_{5}(p)&=&T_{1}(p^2)T_{2}(p^2)+\frac{\alpha_{-2}(p^2)^{2}}{p^2}\nonumber\\
&=&T_{1}(p^2)\bigg(2-\frac{1}{T_{1}(p^2)}\bigg)+\frac{\alpha_{-2}(p^2)^{2}}{p^2}\nonumber\\
&=&2T_{1}(p^2)-1+\frac{\alpha_{-2}(p^2)^{2}}{p^2}\nonumber\\
&=&2\bigg(1+\frac{\alpha_{-2}(p)\delta(p)}{p}+\frac{\alpha_{-2}(p^2)\delta(p^2)}{p^2}\bigg)-1+\frac{\alpha_{-2}(p^2)^{2}}{p^2}\nonumber\\
&=&1-\frac{4\delta(p)}{p}+\frac{2\alpha_{-2}(p^2)\delta(p^2)+\alpha_{-2}(p^2)^{2}}{p^2}.
\end{eqnarray*}
So
\begin{eqnarray*}
&&\!\!\!\!\!\!\!\!\!\!\!\!\!\!\!\!\!\!\!\!\!\!\!\!\!\!\!\!\!T_{1}(p^2)T_{2}(p^2)T_{5}(p)+\frac{\alpha_{-2}(p)^{2}T_{3}(p)T_{4}(p)}{p}=\nonumber\\
&&\!\!\!\!\!1+\frac{4-4\delta(p)}{p}+\frac{\alpha_{-2}(p^2)(\alpha_{-2}(p^2)-4\alpha_{-2}(p^2)\delta(p)+2\delta(p^2))}{p^2}.
\end{eqnarray*}
We note that $\delta(p)=2/(1+1/p)$, $\delta(p^2)=2\delta(p)-1$, and $\alpha_{-2}(p^2)=1$ if $p\leq\sqrt{X}$, $\alpha_{-2}(p^2)=2$ if $\sqrt{X}<p\leq X$. Straightforward calculations then give, if $p\leq\sqrt{X}$,
\begin{equation*}
T_{1}(p^2)T_{2}(p^2)T_{5}(p)+\frac{\alpha_{-2}(p)^{2}T_{3}(p)T_{4}(p)}{p}=\frac{(1-1/p)^{3}}{1+1/p},
\end{equation*}
and, if $\sqrt{X}<p\leq X$,
\begin{equation*}
T_{1}(p^2)T_{2}(p^2)T_{5}(p)+\frac{\alpha_{-2}(p)^{2}T_{3}(p)T_{4}(p)}{p}=\frac{(1-1/p)^{3}}{1+1/p}+O(1/p^2).
\end{equation*}
Hence
\begin{eqnarray*}
S&=&(1+o(1))\prod_{\substack{p\leq\sqrt{X}\\p\nmid q}}\bigg(\frac{(1-1/p)^{3}}{1+1/p}\bigg)\prod_{\substack{\sqrt{X}<p\leq X\\p\nmid q}}\bigg(\frac{(1-1/p)^{3}}{1+1/p}+O(1/p^2)\bigg)\nonumber\\
&=&(1+o(1))\bigg(\prod_{\substack{p\leq X\\p|q}}\frac{(1-1/p)^{3}}{1+1/p}\bigg)^{-1}\prod_{p\leq X}\frac{(1-1/p)^{3}}{1+1/p}\prod_{\sqrt{X}<p\leq X}(1+O(1/p^2))\nonumber\\
&=&(1+o(1))\bigg(\prod_{\substack{p\leq X\\p|q}}\frac{(1-1/p)^{3}}{1+1/p}\bigg)^{-1}\prod_{p\leq X}\bigg(1-\frac{1}{p}\bigg)^{4}\prod_{p}\bigg(1-\frac{1}{p^2}\bigg)^{-1}\nonumber\\
&=&(1+o(1))\frac{\pi^{2}}{6}\prod_{\substack{p>X\\p|q}}\frac{(1-1/p)^{3}}{1+1/p}\bigg(\prod_{p|q}\frac{(1-1/p)^{3}}{1+1/p}\bigg)^{-1}\bigg(\frac{1}{e^{\gamma}\log X}\bigg)^{4}.
\end{eqnarray*}
Combining with \eqref{22} we have
\begin{equation}\label{23}
I_{1}=(1+o(1))\frac{1}{48}\prod_{\substack{p>X\\p|q}}\frac{(1-1/p)^{3}}{1+1/p}\bigg(\frac{\log q}{e^{\gamma}\log X}\bigg)^{4}.
\end{equation}

We now move on to $I_{2}$. We have
\begin{eqnarray*}
I_{2}\ll\frac{1}{\varphi^{*}(q)}\sum_{\substack{m,n\leq (\log q)^{1/4}\\(mn,q)=1}}\frac{|\alpha_{-2}(m)\alpha_{-2}(n)|}{\sqrt{mn}}\sum_{h|q}\varphi(h)\mu(q/h)^{2}\sum_{\substack{ab,cd\leq Z\\acm\equiv \pm bdn(\textrm{mod}\ h)\\acm\ne bdn\\(abcd,h)=1}}\frac{1}{\sqrt{abcd}}.
\end{eqnarray*}
Let
\begin{equation*}
E_{2}(h)=\sum_{\substack{ab,cd\leq Z\\acm\equiv\pm bdn(\textrm{mod}\ h)\\acm\ne bdn\\(abcd,h)=1}}\frac{1}{\sqrt{abcd}}.
\end{equation*}
We divide the terms $ab,cd\leq Z$ into dyadic blocks. Consider the block $Z_{1}\leq ab<2Z_{1}$ and $Z_{2}\leq cd<2Z_{2}$. If $Z_{1}Z_{2}>h^{19/10}$, then by Lemma 11, the sum over this block is
\begin{equation*}
\ll\frac{1}{\sqrt{Z_{1}Z_{2}}}\frac{Z_{1}Z_{2}mn}{h}(\log Z_{1}Z_{2})^{3}\ll\frac{\sqrt{Z_{1}Z_{2}}mn}{h}(\log q)^3,
\end{equation*}
and is $\ll(Z_{1}Z_{2})^{1/2+\epsilon}(mn)^{1+\epsilon}/h$ if $Z_{1}Z_{2}\leq h^{19/10}$. Summing over all the dyadic blocks we have
\begin{equation*}
E_{2}(h)\ll\frac{Zmn}{h}(\log q)^{3}+(mn)^{1+\epsilon}h^{-1/20+\epsilon}.
\end{equation*}
So
\begin{equation*}
\frac{1}{\varphi^{*}(q)}\sum_{h|q}\varphi(h)\mu(q/h)^{2}E(h)\ll\frac{q(\log q)^3}{\varphi^{*}(q)}mn\ll(\log q)^{3+\epsilon}mn.
\end{equation*}
Thus
\begin{eqnarray}\label{24}
I_{2}&\ll&(\log q)^{3+\epsilon}\sum_{m,n\leq (\log q)^{1/4}}|\alpha_{-2}(m)\alpha_{-2}(n)|\sqrt{mn}\nonumber\\
&\ll&(\log q)^{3+\epsilon}\bigg(\sum_{m\leq(\log q)^{1/4}}m^{1/2+\epsilon}\bigg)^2
\ll(\log q)^{15/4+\epsilon}.
\end{eqnarray}
This and \eqref{23} prove Proposition 3.
\end{proof}

\begin{prop}
Let $\delta>0$. Suppose we have $X,q\rightarrow\infty$ with $X\ll(\log\log q)^{2-\delta}$. Then
\begin{equation*}
J=\frac{1}{\varphi^{*}(q)}\sum_{\chi\ (\emph{mod}\ q)}{\!\!\!\!\!\!\!}^{\textstyle{*}}\ \ \bigg|C(\chi)\sum_{\substack{n\in S(X)\\n\leq (\log q)^{1/4}}}\frac{\alpha_{-2}(n)\chi(n)}{\sqrt{n}}\bigg|^2=o\bigg(\frac{\log q}{\log X}\bigg)^{4}.
\end{equation*}
\end{prop}
\begin{proof}
We have
\begin{eqnarray*}
J&\ll&\frac{1}{\varphi^{*}(q)}\sum_{\chi\ (\textrm{mod}\ q)}\bigg|C(\chi)\sum_{\substack{n\in S(X)\\n\leq (\log q)^{1/4}}}\frac{\alpha_{-2}(n)\chi(n)}{\sqrt{n}}\bigg|^2\nonumber\\
&\ll&\frac{\varphi(q)}{\varphi^{*}(q)}\sum_{\substack{ab,cd>Z\\mn\in S(X)\\m,n\leq (\log q)^{1/4}\\acm\equiv\pm bdn(\textrm{mod}\ q)\\(abcdmn,q)=1}}\frac{|\alpha_{-2}(m)\alpha_{-2}(n)|}{\sqrt{abcdmn}}\sum_{\mathfrak{a}\in\{0,1\}}\bigg|W_{\mathfrak{a}}\bigg(\frac{\pi ab}{q}\bigg)W_{\mathfrak{a}}\bigg(\frac{\pi cd}{q}\bigg)\bigg|.
\end{eqnarray*}
We proceed as in Proposition 3. Let us write the last expression as $J_{1}+J_{2}$, where $J_{1}$ consists of the $acm=bdn$ terms and $J_{2}$ consists of the remaining terms. We first estimate $J_{1}$. In the same way as we dealt with $I_{1}$, we write $m=ugh$, $n=uij$, $a=vjk$, $b=vgl$, $c=wil$ and $d=whk$, where $(gh,ij)=(k,l)=(k,gi)=(l,hj)=1$. Also let $f=kl$. Then
\begin{eqnarray*}
J_{1}&\ll&\frac{\varphi(q)}{\varphi^{*}(q)}\sum_{\substack{ughij\in S(X)\\ugh,uij\leq (\log q)^{1/4}\\(ughij,q)=1\\(gh,ij)=1}}\frac{|\alpha_{-2}(ugh)\alpha_{-2}(uij)|}{ughij}\sum_{(f,q)=1}\frac{2^{\omega(f)}}{f}\nonumber\\
&&\quad\sum_{\mathfrak{a}\in\{0,1\}}\bigg|\sum_{\substack{v>\sqrt{Z/gjf}\\(v,q)=1}}\frac{1}{v}W_{\mathfrak{a}}\bigg(\frac{\pi v^{2}gjf}{q}\bigg)\bigg|\bigg|\sum_{\substack{w>\sqrt{Z/hif}\\(w,q)=1}}\frac{1}{w}W_{\mathfrak{a}}\bigg(\frac{\pi w^{2}hif}{q}\bigg)\bigg|
\end{eqnarray*}
We first evaluate the sum over $f$. Let us consider the case $f>q$. The contribution of these terms to the sum is
\begin{eqnarray*}
&\ll&\sum_{f>q}\frac{2^{\omega(f)}}{f}\bigg(\sum_{\substack{v>\sqrt{Z/gjf}\\(v,q)=1}}\frac{1}{v}\bigg(\frac{ v^{2}gjf}{q}\bigg)^{-2}\bigg)\bigg(\sum_{\substack{w>\sqrt{Z/hif}\\(w,q)=1}}\frac{1}{w}\bigg(\frac{w^{2}hif}{q}\bigg)^{-2}\bigg)\nonumber\\
&\ll&\sum_{f>q}\frac{2^{\omega(f)}}{f}\frac{q^4}{f^4}\ll\log q.
\end{eqnarray*}
For the case $f<q$, the expression in the first bracket of the sum is
\begin{equation*}
\ll\sum_{\substack{\sqrt{q/gjf}>v>\sqrt{Z/gjf}\\(v,q)=1}}\frac{1}{v}+1.
\end{equation*}
This is
\begin{eqnarray*}
\ll\frac{\varphi(q)}{q}(\log(\sqrt{q/gjf})-\log(\sqrt{Z/gjf})+O(1+\log\omega(q)))\ll\frac{\varphi(q)}{q}\omega(q),
\end{eqnarray*}
if $gjf<Z_{0}$, and if $gjf>Z_{0}$, it is
\begin{equation*}
\ll\sum_{\substack{v<5^{\omega(q)}\\(v,q)=1}}\frac{1}{v}+1\ll\frac{\varphi(q)}{q}\omega(q).
\end{equation*}
Thus the contribution of the terms $f<q$ to the sum is
\begin{equation*}
\ll\bigg(\frac{\varphi(q)}{q}\omega(q)\bigg)^{2}\sum_{\substack{f\leq q\\(f,q)=1}}\frac{2^{\omega(f)}}{f}\ll\frac{\varphi(q)}{\varphi^{*}(q)}\bigg(\frac{\varphi(q)}{q}w(q)\bigg)^2\frac{\varphi(q)^2}{q^2}(\log q)^2\ll(\log q)^{2+\epsilon}.
\end{equation*}
Hence
\begin{equation*}
J_{1}\ll\frac{\varphi(q)}{\varphi^{*}(q)}(\log q)^{2+\epsilon}\sum_{\substack{ughij\in S(X)\\ugh,uij\leq (\log q)^{1/4}\\(ughij,q)=1\\(gh,ij)=1}}\frac{|\alpha_{-2}(ugh)\alpha_{-2}(uij)|}{ughij}.
\end{equation*}
By the same method used in estimating $S$ in Proposition 3, the inner sum is $\ll(\log X)^{8}$. Thus
\begin{equation}\label{25}
J_{1}\ll(\log q)^{2+\epsilon}(\log X)^{8}=o\bigg(\frac{\log q}{\log X}\bigg)^4.
\end{equation}

To estimate $J_{2}$, first note that
\begin{eqnarray*}
J_{2}\ll\frac{\varphi(q)}{\varphi^{*}(q)}\sum_{\substack{mn\in S(X)\\m,n\leq (\log q)^{1/4}\\(mn,q)=1}}\frac{|\alpha_{-2}(m)\alpha_{-2}(n)|}{\sqrt{mn}}\sum_{\substack{ab,cd>Z\\acm\equiv\pm bdn(\textrm{mod}\ q)\\acm\ne bdn\\(abcd,q)=1}}\frac{1}{\sqrt{abcd}}\bigg(1+\frac{ab}{q}\bigg)^{-2}\bigg(1+\frac{cd}{q}\bigg)^{-2}.
\end{eqnarray*}
We divide the innermost sum into dyadic blocks $Z_{1}<ab\leq2Z_{1}$, $Z_{2}<cd\leq2Z_{2}$, where $Z_{1},Z_{2}>Z$. We have
\begin{eqnarray*}
&&\sum_{\substack{Z_{1}<ab\leq2Z_{1}\\Z_{2}<cd\leq2Z_{2}\\acm\equiv\pm bdn(\textrm{mod}\ q)\\acm\ne bdn\\(abcd,q)=1}}\frac{1}{\sqrt{abcd}}\bigg(1+\frac{ab}{q}\bigg)^{-2}\bigg(1+\frac{cd}{q}\bigg)^{-2}\ll\nonumber\\
&&\quad\qquad\qquad\ll\frac{1}{\sqrt{Z_{1}Z_{2}}}\bigg(1+\frac{Z_{1}}{q}\bigg)^{-2}\bigg(1+\frac{Z_{2}}{q}\bigg)^{-2}\sum_{\substack{Z_{1}<ab\leq2Z_{1}\\Z_{2}<cd\leq2Z_{2}\\acm\equiv\pm bdn(\textrm{mod}\ q)\\acm\ne bdn\\(abcd,q)=1}}1\nonumber\\
&&\quad\qquad\qquad\ll\bigg(1+\frac{Z_{1}}{q}\bigg)^{-2}\bigg(1+\frac{Z_{2}}{q}\bigg)^{-2}\frac{\sqrt{Z_{1}Z_{2}}mn}{q}(\log Z_{1}Z_{2})^3.
\end{eqnarray*}
Summing over all such blocks we find that the innermost sum in the formula for $J_{2}$ is $\ll(\log q)^{3}mn$. So
\begin{equation}\label{26}
J_{2}\ll\frac{\varphi(q)}{\varphi^{*}(q)}(\log q)^{3}\sum_{m,n\leq(\log q)^{1/4}}|\alpha_{-2}(m)\alpha_{-2}(n)|\sqrt{mn}\ll(\log q)^{15/4+\epsilon}.
\end{equation}
This and \eqref{25} prove Proposition 4.
\end{proof}

The first part of Theorem 4 now easily follows from Proposition 3, Proposition 4, \eqref{20}, Lemma 3 and Cauchy's inequality.

The second part then follows by Lemma 2.

\specialsection*{\textbf{Acknowledgments}}
We are grateful to Tim Browning, Brian Conrey and Steve Gonek for helpful comments and suggestions, and to Zeev Rudnick for a stimulating discussion. J.P.K. is supported by an EPSRC Senior Research Fellowship.

\end{document}